\documentclass[10pt,a4paper]{amsart}

\usepackage{amssymb,amsmath,graphicx,enumerate,xcolor,mathtools}
\usepackage{hyperref} % prolinkuje odkazy na clanky
\hypersetup{hidelinks}
\usepackage[shortlabels]{enumitem}
\numberwithin{equation}{section}

\newtheorem{thm}{Theorem}[section]
\newtheorem{question}[thm]{Question}

\newtheorem{lemma}[thm]{Lemma}
\newtheorem{prop}[thm]{Proposition}
\newtheorem{cor}[thm]{Corollary}

\newtheorem*{claim}{Claim}

\theoremstyle{definition}

\newtheorem{remark}[thm]{Remark}

\newtheorem{thmx}{Theorem}
\newtheorem{corx}{Corollary}

\def\C{\mathcal C}

\let\PolishL\L
\def\L{\mathcal L}
\def\U{\mathcal U}

\def\Id{\operatorname{Id}}

\def\en{\mathbb N}

\def\er{\mathbb R}

\def\spt{\operatorname{supp}}
\DeclareMathOperator{\height}{ht}
\DeclareMathOperator{\Lip}{Lip}

\def\ep{\varepsilon}

\newcommand{\norm}[1]{\left\|#1\right\|}
\newcommand{\abs}[1]{\left| #1  \right|}

\newcommand{\setsep}{:\,}

\begin{document}

\title[The classification by positive isomorphisms]{The classification of $C(K)$ spaces for countable compacta by positive isomorphisms}

\author[M. C\' uth]{Marek C\'uth}
\author[J. Havelka]{Jon\'a\v{s} Havelka}
\author[J. Rondo\v{s}]{Jakub Rondo\v{s}}
\author[B. Sari]{B\"unyamin Sar\i}
\email{marek.cuth@matfyz.cuni.cz}
\email{jonas.havelka@volny.cz}
\email{jakub.rondos@gmail.com}
\email{bunyamin@unt.edu}

\thanks{Corresponding author: Jakub Rondo\v{s};
e-mail: \href{mailto:jakub.rondos@gmail.com}{jakub.rondos@gmail.com}.}

\address[M.~C\' uth, J.~Havelka]{Charles University, Faculty of Mathematics and Physics, Department of Mathematical Analysis, Sokolovsk\'a 83, 186 75 Prague 8, Czech Republic}
\address[J. Rondo\v{s}]{Department of Mathematics, Faculty of Electrical Engineering,
Czech Technical University in Prague, Technick\'a 2, 166 27 Praha 6, Czech
Republic}
\address[B. Sari]{Department of Mathematics, University of North Texas, 1155 Union Circle \#311430, Denton, TX 76203-5017}


%\thanks{M. C\'uth was supported by the Czech Science Foundation, project no. GA\v{C}R 24-10705S.}

\begin{abstract}
    We study the classification of spaces of continuous functions $\C(K)$ under positive linear maps. For infinite countable compacta, we show that whenever $\C(K)$ and $\C(L)$ are isomorphic, there exists an isomorphism $T:\C(K)\to \C(L)$ satisfying either $T\geq 0$ or $T^{-1}\geq 0$. We also prove that for any compact spaces $K$ and $L$, the existence of a positive embedding $T: \C(K) \to \C(L)$ implies that the Cantor--Bendixson height of $K$ does not exceed the height of $L$. Further, we introduce a one-sided positive Banach--Mazur distance and compute it in several families of the corresponding spaces of continuous functions. Our estimates also yield new exact values for the classical Banach–Mazur distance between such spaces.

\end{abstract}

\maketitle

%\tableofcontents

\section*{Introduction}

The study of the isometric and isomorphic properties of $\C(K)$ spaces is an active and well-developed topic within the theory of Banach spaces. Concerning the isometric theory, by the Banach--Stone theorem, $\C(K)$ and $\C(L)$ are linearly isometric if and only if the compact spaces $K$ and $L$ are homeomorphic. Moreover, by the Mazurkiewicz--Sierpi\'nski theorem, for every infinite countable compact space $K$ there is a unique countable ordinal $\alpha$ and unique natural number $n$ such that $K$ is homeomorphic to $[0,\omega^\alpha\cdot n]$. In order to shorten the notation, in what follows given an ordinal $\alpha$ we write $\C(\alpha)$ instead of $\C([0,\alpha])$. Concerning the isomorphic theory, in 1960 Bessaga and Pe{\l}czy\'nski \cite{BessagaPelcynski_classification} proved that for any infinite countable compact space $K$ there exists a unique ordinal $\xi<\omega_1$ such that $\C(K)$ is isomorphic to
$\C(\omega^{\omega^\xi})$ and in 1966 Miljutin \cite{Mil66} proved that whenever $K$ is an uncountable compact metric space, then $\C(K)$ is isomorphic to $\C([0,1])$, which completed isomorphic classification of $\C(K)$ spaces for $K$ metrizable compacta. We refer, for example, to \cite[Chapter 4]{albiacKniha} or \cite[Section 2.6]{HMVZ} for more information and additional historical background. Both qualitative and quantitative variants of the results mentioned above were considered by many researchers, especially during the last 20 years, some of which are mentioned below.

In the present work, we consider the following strengthening of the notion of an isomorphism: we say that a bounded linear operator $T:\C(K)\to \C(L)$ is \emph{positive embedding} (resp. \emph{positive isomorphism}) if $T$ is an isomorphic embedding (resp. surjective isomorphism) with $T\geq 0$ (that is, $Tf\geq 0$ whenever $f\geq 0$). Let us emphasize that we do not assume positive isomorphism $T$ satisfies $T^{-1}\geq 0$. Indeed, by the well-known Kaplansky's theorem \cite{K47}, if both $T\geq 0$ and $T^{-1}\geq 0$, then $K$ and $L$ are homeomorphic, and hence $\C(K)$ and $\C(L)$ are lattice isometric. We shall sometimes write (one-sided) isomorphic embedding and (one-sided) positive isomorphism in order to emphasize that we do not assume $T^{-1}\geq 0$. Our first main result could be understood as a qualitative improvement of the result by Bessaga and Pe{\l}czy\'nski mentioned above. The symbol $\height(K)$ used below denotes the height of the compact space $K$, see Section~\ref{sec:prelim} for a precise definition.

\begin{thmx}\label{thm:Intro1}  Let $K$ and $L$ be infinite countable compact spaces. The following are equivalent:
\begin{enumerate}[label=(\roman*)]
   \item There exists a positive isomorphism from $C(K)$ onto $C(L)$.
    \item $C(K)$ is isomorphic to $C(L)$ and $\height(K)\leq \height(L)$.
\end{enumerate}

Moreover, for any compact spaces $K$ and $L$, if there exists a positive embedding $T:C(K)\to C(L)$, then $\height(K)\leq \height(L)$.
\end{thmx}

A notable consequence of Theorem~\ref{thm:Intro1} is that given infinite countable compacta $K$ and $L$, the Banach spaces $\C(K)$ and $\C(L)$ are isomorphic if and only if there exists a surjective isomorphism $T:\C(K)\to \C(L)$ such that either $T\geq 0$ or $T^{-1}\geq 0$. A second notable consequence is the following (the proof follows from Theorem~\ref{thm:Intro1} and certain earlier results, see Section~\ref{sec:main} for more details).

\begin{corx}
\label{cor:main1}
Let $K$ and $L$ be compact spaces, and $K$ is metrizable. Then there exists a positive embedding of $\C(K)$ into $\C(L)$ if and only if $ht(K) \leq ht(L)$.   
\end{corx}

Moreover, some parts of our argument produce quantitative estimates that seem to be new even in the classical context. Recall that given compact spaces $K$ and $L$ we define the \emph{Banach--Mazur distance} between $\C(K)$ and $\C(L)$, denoted $d_{BM}(\C(K),\C(L))$, as the infimum of the numbers $\norm{T}\cdot \norm{T^{-1}}$ over all surjective isomorphisms $T:\C(K)\to \C(L)$. Similarly, we define $d_{BM}^+(\C(K),\C(L))$ to be the infimum of the numbers $\norm{T}\cdot \norm{T^{-1}}$ over all positive isomorphisms $T:\C(K)\to \C(L)$. Notice that $d_{BM}^+$ is not symmetric: by Theorem~\ref{thm:Intro1}, for example, $d_{BM}^+(\C(\omega),\C(\omega^2))<\infty$ while $d_{BM}^+(\C(\omega^2),\C(\omega))=\infty$. Once we obtain an estimate of $d_{BM}$ or $d_{BM}^+$, it is natural to ask whether this estimate is optimal. The following summarizes some of the most interesting quantitative estimates we obtain.

\begin{thmx}\label{thm:Intro2}
Let $\alpha$ be a countable nonzero ordinal and $n\in\en$.
\begin{enumerate}[label=(\roman*)]
\item We have
\[
d_{BM}(\C(\omega^{\alpha}), \C(\omega^{\alpha}\cdot n))\leq d_{BM}^+(\C(\omega^{\alpha}), \C(\omega^{\alpha}\cdot n)) \leq 2+\sqrt{5}.
\]
\item\label{it:itemMainPower}
If $\alpha = \omega^\beta$ for some ordinal $\beta$, then
\[
d_{BM}(\C(\omega^{\alpha}), \C(\omega^{\alpha\cdot n})) = d_{BM}^+(\C(\omega^{\alpha}), \C(\omega^{\alpha\cdot n})) = n+\sqrt{(n-1)(n+3)}. 
\]
\end{enumerate}
\end{thmx}

\begin{thmx}\label{thm:Intro3}
For every nonzero countable ordinal $\alpha$ and every integer $k\geq2$,
\[
d_{BM}^+(\C(\omega^\alpha\cdot k),\C(\omega^\alpha))
=k+\sqrt{(k-1)(k+3)}.
\]
In the opposite orientation, if $\alpha=\omega^\beta$ for some countable
ordinal $\beta$, then
\[
d_{BM}^+(\C(\omega^\alpha),\C(\omega^\alpha\cdot2))=2+\sqrt3,
\]
and, for every $k\geq3$,
\[
d_{BM}^+(\C(\omega^\alpha),\C(\omega^\alpha\cdot k))=2+\sqrt5.
\]
\end{thmx}

Let us now discuss the relevance of Theorems~\ref{thm:Intro1}--\ref{thm:Intro3} to some recent results and formulate open problems we were not able to solve. There are several existing works in which the authors consider various morphisms between $\C(K)$ spaces, with the primary interest being improvements of the Banach--Stone theorem. By Milgram's theorem \cite{M49} (resp. Kaplansky's theorem \cite{K47}, resp. Jarosz's theorem \cite{J90}), if there exists a bijection $T:\C(K)\to \C(L)$ such that both $T$ and $T^{-1}$ preserve multiplication (resp. pointwise ordering, resp. disjointness of supports), then $K$ and $L$ are homeomorphic. There is a recent attempt to find a common generalization of those results by Kania and Rmoutil using the notion of compatibility isomorphism, see \cite{KR18} and \cite{KR18errata}. Another direction in generalizing the Banach--Stone theorem is provided by results of Dutrieux and Kalton \cite{DK05}, and of Galego and da Silva \cite{GS16}, who showed that if the distance (Gromov--Hausdorff, Kadets, Lipschitz, or even coarse) between $\C(K)$ spaces is small, then the corresponding Banach spaces are linearly isometric. Positive isometric embeddings between general $\C(K)$ spaces were studied also in \cite{Pl13}. Our Theorem~\ref{thm:Intro1} shows that a one-sided variant of Kaplansky's theorem yields an isomorphic characterization of $\C(K)$ spaces when $K$ is an infinite countable compact space. This suggests further directions of research, including one-sided variants of some of the other results mentioned above. The most pressing open problem for us is whether Theorem~\ref{thm:Intro1} extends to uncountable compact metric spaces, that is, whether Miljutin's theorem can be improved in the following sense.

\begin{question}\label{qe:milj}
    Let $K$ be an uncountable compact metric space. Is it true that there exists an isomorphism $T:\C(K)\to \C([0,1])$ such that either $T\geq 0$ or $T^{-1}\geq 0$?
\end{question}

Theorem~\ref{thm:Intro2} is very much related to the recent efforts of finding the Banach--Mazur distance between various $\C(K)$ spaces. This is motivated by the problem of Bessaga and Pe{\l}czy\'nski \cite{BessagaPelcynski_classification},  posed in 1960, in which they asked for the exact value of the Banach--Mazur distance between $\C(\alpha)$ and $\C(\beta)$ for ordinals $\alpha,\beta$.  Pe{\l}czy\'nski even raised a conjecture in 1968 \cite[Problem 28]{Pel68} as to whether this distance must always be an integer. In 1970 Gordon \cite{Gordon3} proved that the distance between non-isometric $\C(K)$ spaces with $K$ countable has to be at least $3$ and $d_{BM}(\C(\omega),\C(\omega\cdot 2)) = 3$. Since then, the problem remained basically untouched until 2013, when Galego with Candido \cite{CandidoGalegoComega} proved (among other results) that $d_{BM}(\C(\omega),\C(\omega\cdot n))\leq 2+\sqrt{5}$ and $d_{BM}(\C(\omega),\C(\omega^n))\leq n + \sqrt{(n-1)(n+3)}$ for every $n\geq 2$. In 2025, Malec and Piasecki \cite{MaPi25} proved that equality actually holds for the latter, and we have $d_{BM}(\C(\omega),\C(\omega^n)) = n + \sqrt{(n-1)(n+3)}$. Our Theorem~\ref{thm:Intro2} improves and extends those results (we note that proofs required a number of new techniques and ideas and we refer the interested reader to Remark~\ref{rem:novelBesPel}, Remark~\ref{rem:upperBound} and Remark~\ref{rem:piasAnalogy} for more details).

We should mention that the issue of estimating the value of $d_{BM}(\C(\alpha),\C(\beta))$ for various ordinals $\alpha,\beta$ was undertaken also in the recent paper by Somaglia and the last author of this paper, see \cite{rondos-somaglia}. We find it quite interesting that despite the efforts of several researchers the problem of finding the exact value of $D_3:=d_{BM}(\C(\omega),\C(\omega\cdot 3))$ remains open. In 2024 Gergont and Piasecki \cite{Piasecki_distance_not_integer} proved $D_3\in (3.23,3.88)$, disproving the conjecture by Pe{\l}czy\'nski from 1968 mentioned above, and in a very recent preprint Korpalski and Plebanek \cite{KP25} found a sharper estimate, which gives a better lower bound.

An earlier version of this paper posed the following question. It is answered completely by Theorem~\ref{thm:Intro3}.

\begin{question}\label{que:exact-positive-omega}
    What are the values of $d_{BM}^+(\C(\omega),\C(\omega\cdot 2))$ and $d_{BM}^+(\C(\omega\cdot 2),\C(\omega))$?
\end{question}

The same earlier version also asked whether the classical and positive distances can be substantially different. We retain the question in order to record its role in the subsequent development.

\begin{question}\label{que:positive-separation}
    Are there countable ordinals $\omega\leq \alpha < \beta$ such that $d_{BM}(\C(\alpha),\C(\beta)) < \min\{d_{BM}^+(\C(\alpha),\C(\beta)),d_{BM}^+(\C(\beta),\C(\alpha))\}$?
\end{question}

After the first version of this paper appeared, Korpalski and Plebanek \cite{KPpositive26} were the first to answer Question~\ref{que:positive-separation} affirmatively; they proved that both one-sided positive distances for the pair $\C(\omega)$ and $\C(\omega\cdot2)$ exceed the classical Banach--Mazur distance $3$ \cite{KPpositive26}. Theorem~\ref{thm:Intro3} strengthens their result by determining the two distances exactly. In particular,
\[
3=d_{BM}(\C(\omega),\C(\omega\cdot2))
<2+\sqrt3
=d_{BM}^+(\C(\omega),\C(\omega\cdot2))<2+\sqrt{5}=
d_{BM}^+(\C(\omega\cdot2),\C(\omega)).
\]

The structure of our paper is as follows. In Section~\ref{sec:prelim} we gather our notation and preliminary results. In Section~\ref{sec:main} we prove Theorems~\ref{thm:embedMain} and \ref{thm:class}, from which Theorem~\ref{thm:Intro1} follows together with Corollary~\ref{cor:main1} and several quantitative estimates. In Section~\ref{sec:estimates} we obtain further estimates for $d_{BM}$ and $d_{BM}^+$, including Theorem~\ref{thm:Intro2}, and collect the upper constructions needed later. Finally, Section~\ref{sec:positive-distances} determines the exact positive distances stated in Theorem~\ref{thm:Intro3}.

\section{Preliminaries}\label{sec:prelim}

All topological spaces are tacitly assumed to be Hausdorff. If $K$ is a topological space and $\alpha$ is an ordinal, by $K^{(\alpha)}$ we denote the standard \emph{Cantor--Bendixson derivative of $K$ of order $\alpha$}, that is, we have $K^{(0)}=K$, 
\[K^{(1)}=\{x \in K\setsep x \text{ is an accumulation point of } K\},\]
and for each ordinal $\alpha$, $K^{(\alpha+1)}=(K^{(\alpha)})^{(1)}$, and $K^{(\alpha)}=\bigcap_{\beta < \alpha} K^{(\beta)}$, if $\alpha$ is a limit ordinal. We note that it is easy to check that if for each open set $U \subseteq K$ and an ordinal $\alpha$, it holds $U^{(\alpha)}=U \cap K^{(\alpha)}$, when $U$ is considered as a topological space with the inherited topology, a fact that we will quite frequently use. 

A topological space $K$ is called \emph{scattered} if $K^{(\alpha)}=\emptyset$ for some ordinal $\alpha$. In this case its \emph{height}, denoted by $\height(K)$, is the least such ordinal; for a nonscattered space we put $\height(K)=\infty$, with every ordinal smaller than $\infty$. By compactness, the height of every scattered compact space is a successor ordinal. The existence of a homeomorphism of compact spaces $K, L$ is denoted by $K \simeq L$.

We also recall that each countable compact space is scattered, and moreover, the homeomorphism classes of such spaces are determined by their heights and the cardinality of the last nonempty derivative, by the Mazurkiewicz--Sierpi\'nski classification (see e.g. \cite[Theorem 8.6.10]{semadeni}). Also, it is a standard fact, which can be easily proved by transfinite induction, that for each ordinal $\alpha$ and $n \in \en$, denoting $K=[1, \omega^{\alpha}\cdot n]$, it holds
\[K^{(\alpha)}=\{\omega^{\alpha}i\setsep i=1, \ldots, n\},\]
hence $\height(K)=\alpha+1$ and $\abs{K^{(\height(K)-1)}}=n$. Together with the aforementioned classification, it follows that each countable compact space is homeomorphic to the ordinal interval $[1, \omega^{\alpha}\cdot n]$ for some countable ordinal $\alpha$ and natural number $n \in \en$.

Next, if $K$ is compact, $\C(K)$ stands as usual for the Banach space of continuous real-valued functions on $K$ endowed with the supremum norm.

\section{The classification}\label{sec:main}

The main purpose of this section is to prove our main classification results. 

\subsection{Preservation of height}\label{subsec:height}
We begin by establishing the fundamental rigidity of the topological height under positive embeddings. 

% --- The Main Theorem (formerly Theorem 2.1 / Prop 2.9c) ---
\begin{thm}\label{thm:embedMain}
    Let $K$ and $L$ be compact spaces. If there exists a positive embedding of $\C(K)$ into $\C(L)$, then $\height(K) \leq \height(L)$.
\end{thm}

We start with a few auxiliary results.
Given a compact space $K$ and a family $\U$ of nonempty subsets of $K$, the \emph{order} of $\U$ is the smallest integer $k$ (if it exists) such that each point of $K$ lies in at most $k$ sets from $\U$. The family $\U$ is \emph{point-bounded} if it has order $k$ for some $k \in \en$.

\begin{lemma}
\label{point-bounded functions}
Let $K, L$ be compact spaces and $T:\C(K) \to \C(L)$ be a positive embedding. Assume that $\mathcal{U}$ is a point-bounded system of open subsets of $K$, and for each $U \in \mathcal{U}$, let $f_U \in \C(K, [0, 1])$ be a norm-$1$ function which vanishes outside of the set $U$. Then for each $0<\ep < \frac{1}{\norm{T^{-1}}}$,  the family $\{y \in L\setsep Tf_U(y) \geq \ep\}_{U \in \mathcal{U}}$
is a point-bounded system of compact sets that have nonempty interiors.
%In fact, if $n$ is the order of the system $\mathcal{U}$, then the order of $\{y \in L\setsep Tf_U(y) \geq \ep\}_{U \in \mathcal{U}}$ is at most $\frac{\norm{T}n}{\ep}$.
\end{lemma}

\begin{proof}
The sets forming the system are clearly compact. To show that they have nonempty interiors it is enough to show that for each $U \in \mathcal{U}$ and for each $0<\ep < \frac{1}{\norm{T^{-1}}}$, the set $\{y \in L: Tf_U(y)>\ep \}$ is nonempty. 
This however follows from the elementary estimates
\[\ep<\frac{1}{\norm{T^{-1}}} \leq \norm{Tf_U}=\sup \{\abs{Tf_U(y)}\setsep y \in L\}=\sup \{Tf_U(y)\setsep y \in L\},\]
since $f_U$ is a positive function of norm $1$ and the embedding $T$ is positive.

To proceed further, let $n$ denote the order of the system $\mathcal{U}$ and assume that $k \in \en$ is such that there exists a point $y \in L$ and sets $U_1, \ldots, U_k$ from $\mathcal{U}$ such that $Tf_{U_i}(y) \geq \ep$ for each $i=1, \ldots, k$. Then, we have 

\[k\ep \leq \sum_{i=1}^k Tf_{U_i}(y) \leq \norm{T}\norm{\sum_{i=1}^k f_{U_i}} \leq \norm{T}n .\]
Hence $k \leq \frac{\norm{T}n}{\ep}$, as desired.
\end{proof}

Now, we collect several general properties of positive embeddings of $\C(K)$ spaces related to the Cantor--Bendixson derivatives of $K$, which are also used in Section~\ref{sec:estimates}.

\begin{prop}
\label{prop:basic properties of positive embeddings}
Let $K$ and $L$ be compact spaces and $T:\C(K) \to \C(L)$ be a positive embedding.
    \begin{enumerate}[label=(\alph*)]
    
    \item\label{it:tech_small_image} Assume $U \subseteq K$ is open and $S \subseteq L$ is compact with $\height(S) < \height(U)$. Then for each $0 < \ep < \frac{1}{\norm{T^{-1}}}$, there exists a norm-$1$ function $f\in\C(K,[0,1])$ supported on $U$ such that $\norm{Tf|_S} < \ep$.
        \item\label{it:tech_preservation} For each $0 < \ep < \frac{1}{\norm{T^{-1}}}$ and any function $f \in \C(K, [0, 1])$, if the set $\{x \in K : f(x)=1\}$ contains an open set $U$ with $U^{(\alpha)} \neq \emptyset$, then
        \[ \{y \in L : Tf(y) \geq \ep\}^{(\alpha)} \neq \emptyset. \]
        
    \end{enumerate}
\end{prop}

\begin{proof}
We first prove \ref{it:tech_small_image} by contradiction. Suppose that for every norm-$1$ function $f \in \C(K, [0, 1])$ supported on $U$, we have $\norm{Tf_{|S}} \geq \ep$.
We claim that for every ordinal $\alpha$, if a function $g \in \C(K, [0,1])$ is supported on $U$ and satisfies $g=1$ on an open set $V \subseteq U$ with $V^{(\alpha)} \neq \emptyset$, then
\[ \{y \in S \setsep Tg(y) \geq \ep\}^{(\alpha)} \neq \emptyset. \]
If this claim holds, then taking $\alpha = \height(S)$ yields a contradiction, as $S^{(\height(S))} = \emptyset$ by definition.

We proceed by transfinite induction on $\alpha$. The base case $\alpha=0$ follows immediately from the assumption that $\norm{Tg_{|S}} \geq \ep$.
Assume the claim holds for $\alpha$. Let $g$ be supported on $U$ with $g=1$ on $V$, where $V^{(\alpha+1)} \neq \emptyset$.
Choose a sequence of distinct points $(x_n)_{n \in \en} \subseteq V^{(\alpha)}$ and a corresponding sequence of pairwise disjoint open neighborhoods $(U_n)_{n \in \en} \subseteq V$ such that $x_n \in U_n$.
For each $n$, let $g_n \in \C(K, [0, 1])$ be a function supported on $U_n$ such that $g_n=1$ on a neighborhood of $x_n$.
By the inductive hypothesis, for each $n$ there exists a point $y_n \in \{y \in S \setsep Tg_n(y) \geq \ep\}^{(\alpha)}$.

By Lemma \ref{point-bounded functions}, the collection of sets $\{ \{y \in S \setsep Tg_n(y) \geq \ep\} \}_{n \in \en}$ is point-bounded. Consequently, the collection of their $\alpha$-derivatives is also point-bounded.
Since the points $y_n$ belong to at most finitely many of these sets, the set $\bigcup_{n=1}^\infty \{y \in S \setsep Tg_n(y) \geq \ep\}^{(\alpha)}$ must be infinite.
Using the positivity of $T$ and $g_n \leq g$, we obtain $Tg_n \leq Tg$. Therefore,
\[ \{y \in S \setsep Tg(y) \geq \ep\}^{(\alpha)} \supseteq \bigcup_{n=1}^\infty \{y \in S \setsep Tg_n(y) \geq \ep\}^{(\alpha)}. \]
Since the set on the left is compact and contains an infinite subset, it must have a nonempty derivative. Thus, $\{y \in S \setsep Tg(y) \geq \ep\}^{(\alpha+1)} \neq \emptyset$. Hence, the claim holds for $\alpha + 1$.

If $\alpha$ is a limit ordinal and the claim holds for each $\beta<\alpha$, then the claim holds for $\alpha$ by compactness using the following intersection of nested compact sets:
\[ \{y \in S \setsep Tg(y) \geq \ep\}^{(\alpha)} = \bigcap_{\beta < \alpha} \{y \in S \setsep Tg(y) \geq \ep\}^{(\beta)} \neq \emptyset. \]
This completes the proof of \ref{it:tech_small_image}.

We now derive \ref{it:tech_preservation} from \ref{it:tech_small_image}. Suppose there exists a function $f$ satisfying the hypothesis, but the set $K_f = \{y \in L : Tf(y) \geq \ep\}$ satisfies $K_f^{(\alpha)} = \emptyset$.
Let $S = K_f$. Then $\height(S) \leq \alpha$.
Since the support of $f$ contains an open set $U$ with $U^{(\alpha)} \neq \emptyset$, we have $\height(U) > \alpha \geq \height(S)$.
We can now apply \ref{it:tech_small_image} to find a norm-$1$ function $h \in \C(K, [0,1])$ supported on $U$ such that $\norm{Th_{|S}} < \ep$.
Since $h \leq f$ and $T$ is positive, we have $Th \leq Tf$.
Consequently, if $y \notin S$, then $Tf(y) < \ep$, and thus $Th(y) < \ep$.
Combined with the fact that $\norm{Th_{|S}} < \ep$, we obtain $\norm{Th} < \ep$.
This contradicts the condition $\ep < \frac{1}{\norm{T^{-1}}}$.
\end{proof}

\begin{proof}[Proof of Theorem~\ref{thm:embedMain}]
    Apply Proposition~\ref{prop:basic properties of positive embeddings} to the constant function $f \equiv \mathbf{1}$ on $K$. 
\end{proof}

\subsection{Positive classification of spaces on countable compacta}

In order to formulate our second step towards the proof of Theorem~\ref{thm:Intro1}, let us recall the following notation. For an ordinal $\alpha$, $\Gamma(\alpha)$ stands for the minimal ordinal of the form $\omega^{\beta}$ which is not less than $\alpha$. That is, for any ordinal $\eta$ we have $\Gamma(\omega^{\eta}) = \omega^\eta$ and $\Gamma(\alpha)=\omega^{\eta+1}$ for ordinal $\alpha$ satisfying $\omega^\eta<\alpha\le\omega^{\eta+1}$.

\begin{thm}
\label{thm:class}
Let $K$ and $L$ be infinite countable compact spaces. Then the following assertions are equivalent:
\begin{enumerate}[label=(\roman*)]
    \item\label{it:positivelyIso} There exists a positive isomorphism from $C(K)$ onto $C(L)$,
    \item\label{it:isoAndHeight} $\C(K)$ is isomorphic to $\C(L)$ and $\height(K) \leq \height(L)$,
    \item\label{it:height} $\height(K) \leq \height(L) < \Gamma(\height(K))$.
\end{enumerate}
\end{thm}

We note that while proving \ref{it:height} implies \ref{it:positivelyIso} in Theorem~\ref{thm:class} we also obtain quantitative estimates for the (one-sided) positive Banach--Mazur distance, some of which are new even in the classical setting. Some of those are mentioned in Theorem~\ref{thm:Intro2} above, for all the estimates we obtain in this section see Proposition~\ref{prop:Construction with better bound}, Proposition~\ref{prop:constructing-isomorphism-new-attempt} and Corollary~\ref{cor:estimatesForOmegaToAlpha}.

Before discussing proofs of our main results, let us mention some of the consequences. The first one is Corollary~\ref{cor:main1} stated in the introduction.

\begin{proof}[Proof of Corollary~\ref{cor:main1}] The implication that if there exists a positive embedding of $\C(K)$ into $\C(L)$, then $\height(K) \leq \height(L)$ is shown in Theorem~\ref{thm:embedMain}. 

For the other implication, the case when $K$ is nonscattered follows from classical results. For the proof as well as related references, see \cite[Theorem 4.1]{Rondos-Sobota_copies_of_separable_C(L)}. If $K$ is scattered, then it is homeomorphic to the ordinal interval $[1, \omega^{\alpha}\cdot n]$ for some countable ordinal $\alpha$ and $n \in \en$. Then, again by \cite[Theorem 4.1]{Rondos-Sobota_copies_of_separable_C(L)}, there exists a positive isometric embedding of $\C([1, \omega^{\alpha}])$ into $\C(L)$. Further, by Theorem~\ref{thm:class}, there exists a positive isomorphism of $\C([1, \omega^{\alpha}\cdot n])$ onto $\C([1, \omega^{\alpha}])$, which finishes the proof.
\end{proof}

From the above proof, we get also the following equivalence (which, in the case of nonscattered $K$, again follows from classical results).

\begin{cor}
Let $K, L$ be compact spaces such that $K$ is metrizable. Suppose that either $K$ is nonscattered, or $K$ is scattered and $\abs{K^{(\height(K)-1)}}=1$. Then, the following assertions are equivalent:
\begin{itemize}
    \item[(i)] $\height(K) \leq \height(L)$,
    \item[(ii)] $\C(K)$ embeds into $\C(L)$ by a positive isomorphism,
    \item[(iii)] $\C(K)$ embeds into $\C(L)$ by a positive isometry.
\end{itemize}
\end{cor}

As an additional consequence of Theorem~\ref{thm:class} we have the following.

\begin{cor}
 \label{main}   
For countable compact spaces $K, L$ the following assertions are equivalent:
\begin{itemize}
    \item[(i)] $\height(K)=\height(L)$,
    \item[(ii)] There is a positive isomorphism from $\C(K)$ onto $\C(L)$ and also a positive isomorphism from $\C(L)$ onto $\C(K)$.
    \item[(iii)] There is a positive embedding from $\C(K)$ into $\C(L)$ and also a positive embedding from $\C(L)$ into $\C(K)$.
    \end{itemize}
\end{cor}
\begin{proof}
(i) $\Rightarrow$ (ii) follows from Theorem~\ref{thm:class}, the implication (ii) $\Rightarrow$ (iii) is trivial, and the implication (iii) $\Rightarrow$ (i) follows from Corollary~\ref{cor:main1}.   
\end{proof}

Now, let us discuss the proof of Theorem~\ref{thm:class}. We start with the following remark.

\begin{remark}\label{rem:szlenkAndIso}
    The equivalence between \ref{it:isoAndHeight} and \ref{it:height} is a direct consequence of the classical classification of $\C(K)$ spaces for countable $K$ by Bessaga and Pe{\l}czy\'nski \cite{BessagaPelcynski_classification}. Recall that for infinite countable ordinals $\alpha$ and $\beta$ with $\alpha \leq \beta$, the spaces $\C(\alpha)$ and $\C(\beta)$ are isomorphic if and only if $\beta < \alpha^\omega$. By basic ordinal arithmetic, this condition is equivalent to $\Gamma(\height(\alpha)) = \Gamma(\height(\beta))$. Consequently, under the assumption $\height(K) \leq \height(L)$, isomorphism holds if and only if $\height(L) < \Gamma(\height(K))$.
\end{remark}

The fact that \ref{it:positivelyIso} implies \ref{it:isoAndHeight} follows immediately from Theorem~\ref{thm:embedMain}. Concerning the implication that \ref{it:height} implies \ref{it:positivelyIso} we note the following.

\begin{remark}\label{rem:novelBesPel} Our proof of the implication \ref{it:height} $\Rightarrow$ \ref{it:positivelyIso} relies on a direct construction based on a decomposition of ordinal intervals, diverging from the classical  method of Bessaga and Pe{\l}czy\'nski (see, e.g., \cite[Lemma 2.57]{HMVZ}). The classical argument utilizes the isomorphism between $\C(K)$ and its $c_0$-sum. However, such a strategy is obstructed in the positive setting because $c$ is not positively isomorphic to $c_0$ (see Proposition~\ref{prop:cAndc0} below). Consequently, a constructive approach appears necessary. Furthermore, this direct construction yields the explicit quantitative estimates for the positive Banach--Mazur distance formulated in Propositions~\ref{prop:Construction with better bound} and \ref{prop:constructing-isomorphism-new-attempt}. We would like to note that we made a considerable effort to obtain optimal constants in Proposition~\ref{prop:Construction with better bound} and Proposition~\ref{prop:constructing-isomorphism-new-attempt}, if we were interested only in qualitative results, the proofs could be substantially simplified.
\end{remark}

\begin{prop}\label{prop:cAndc0}
    There exists no order-preserving bijection $T: c \to c_0$ such that $T^{-1}$ is Lipschitz continuous.
\end{prop}

\begin{proof}
    Suppose for the sake of contradiction that such a map $T$ exists. By replacing $T$ with $f \mapsto T(f) - T(0)$, we may assume without loss of generality that $T(0)=0$. Let $L = \Lip(T^{-1})$.
    Consider the constant sequence $\mathbf{1} \in c$. Since $T(\mathbf{1}) \in c_0$, there exists $k \in \en$ such that 
    \[ (T\mathbf{1})_k < \frac{1}{L}. \]
    Define $x = \frac{1}{L}e_k \in c_0$, where $e_k$ is the standard basis vector. Then $\norm{x} = 1/L$, and the Lipschitz condition implies
    \[ \norm{T^{-1}x} = \norm{T^{-1}x - T^{-1}0} \leq L \norm{x} = 1. \]
    In the lattice structure of $c$, $\norm{y} \leq 1$ implies $y \leq \mathbf{1}$. Therefore, $T^{-1}x \leq \mathbf{1}$. Since $T$ is order-preserving, this implies $x \leq T\mathbf{1}$. However, evaluating at the $k$-th coordinate yields
    \[ x_k = \frac{1}{L} > (T\mathbf{1})_k, \]
    which contradicts the inequality $x \leq T\mathbf{1}$.
\end{proof}

  \subsection{Construction of the positive isomorphisms}

In this subsection, we construct the positive isomorphisms required to establish the implication \ref{it:height} $\Rightarrow$ \ref{it:positivelyIso} of Theorem~\ref{thm:class}. We note that, since homeomorphism classes of compact spaces correspond to Banach-lattice isometries of the respective spaces of continuous functions, we can instead construct positive isomorphisms of spaces of the form $\C(\tilde{K})$, $\C(\tilde{L})$, where $\tilde{K} \simeq K$ and $\tilde{L} \simeq L$, which will enable us to simplify the notation several times. Our construction relies on a specific family of operators whose distortion is governed by a set of weights. We begin by determining the optimal constant associated with these weights, which will serve as a bound for the Banach--Mazur distance in our main results. 

For $n \in \en$, let $\Delta_n = \{(\lambda_1, \dots, \lambda_n) \in (0, 1]^n : \sum_{i=1}^n \lambda_i = 1\}$. We define the constant $C(n)$ as
\[
C(n) = \inf_{\lambda \in \Delta_n} \left\{ \max \left( \frac{2}{\lambda_1}-1, \max_{i=2, \dots, n}\left(\frac{2}{\lambda_i}+1\right) \right) \right\}.
\]

\begin{lemma}
\label{lemma:the constant}
For each $n \in \en$, $C(n)=n+\sqrt{(n-1)(n+3)}$, and the infimum is attained.
\end{lemma}
 \begin{proof}
    Let
    \[
    F(\lambda) \coloneqq \max\left(\frac{2}{\lambda_1}-1,\ \max_{i=2,\dots,n}\left(\frac{2}{\lambda_i}+1\right)\right), \qquad \lambda \in \Delta_n.
    \]

    Fix $\lambda_1$ and let $m \coloneqq \min_{2\le i\le n}\lambda_i$. Since $\sum_{i=2}^n \lambda_i \ge (n-1)m$, we have
    \[
    m \le \frac{1-\lambda_1}{n-1}.
    \]
    The function $t \mapsto \frac{2}{t}+1$ is strictly decreasing on $(0,1]$.  We have
    \[
    \max_{i=2,\dots,n}\left(\frac{2}{\lambda_i}+1\right) = \frac{2}{m}+1 \ge  \frac{2(n-1)}{1-\lambda_1}+1.
    \]
    Moreover, equality holds for $\lambda_2=\cdots=\lambda_n=\frac{1-\lambda_1}{n-1}$. Therefore, 
    \[
    \inf_{\substack{\lambda\in \Delta_n\\\lambda_1=t}}F(\lambda) =\max\left(\frac{2}{t}-1,\ \frac{2(n-1)}{1-t}+1\right).
    \]
    
    Thus $$C(n)=\inf_{t\in (0,1)}\max\left(\frac{2}{t}-1,\ \frac{2(n-1)}{1-t}+1\right).$$
    Since $\frac{2}{t}-1$ is decreasing and $\frac{2(n-1)}{1-t}+1$ is increasing, the quantity above is minimized at the unique point $x \in (0,1)$ where
    \[
    \frac{2}{x}-1 = \frac{2(n-1)}{1-x}+1.
    \]
    This yields $$x^2 - (n+1)x + 1 = 0.$$
    
    The only root in $(0,1)$ is $x=\frac{(n+1)-\sqrt{(n-1)(n+3)}}{2}$. Thus 
    $$C(n)= n+\sqrt{(n-1)(n+3)}.$$
   
   This infimum is achieved by setting $\lambda_1 = x$ and $\lambda_2 = \dots = \lambda_n = \frac{1-x}{n-1}$ where $x$ is the root above.
\end{proof}

\begin{lemma}
\label{lemma:split into subintervals}
Let $\alpha$ be a nonzero countable ordinal and let $\Gamma$ be a nonempty
at most countable set. Then there are pairwise disjoint sets
$\{I_\gamma\setsep \gamma\in\Gamma\}$ such that $\bigcup_{\gamma\in\Gamma} I_\gamma=[0,\omega^\alpha)$, and the following conditions are satisfied.
\begin{enumerate}[label=(i-\alph*)]
\item\label{it:splitIntoIntervals}
For every $\gamma\in\Gamma$, the set
$I_\gamma\subset[0,\omega^\alpha)$ is open, homeomorphic to
$[0,\omega^\alpha)$, and $\overline{I_\gamma}=I_\gamma\cup\{\omega^\alpha\}$.

\item\label{it:separateIn}
For every $\eta<\omega^\alpha$, the set
$\{\gamma\in\Gamma\setsep \min I_\gamma\leq\eta\}$ is finite. In particular, if $(\gamma_n)$ is any one-to-one sequence in $\Gamma$, then
$\lim_{n\to\infty}\min I_{\gamma_n}=\omega^\alpha$.
\end{enumerate}
\end{lemma}
\begin{proof}
We shall define each $I_\gamma$ as a disjoint union of compact intervals.
Choose a sequence $(\sigma_m)_{m\in\en}$ in $\Gamma$ such that every
$\gamma\in\Gamma$ occurs in this sequence infinitely many times.

First, let us consider the case when $\alpha=\beta+1$ for some countable
ordinal $\beta$. Put
\[
   J_1:=[0,\omega^\beta]\quad \text{and}\quad J_m:=
   [\omega^\beta\cdot(m-1)+1,\omega^\beta\cdot m] \text{ for } m\geq 2. 
\]

Now, consider the case when $\alpha$ is a limit ordinal. Pick a strictly
increasing sequence $(\beta_m)$ of ordinals converging to $\alpha$ and put
\[
   J_1:=[0,\omega^{\beta_1}]\quad \text{and}\quad J_m:=
   [\omega^{\beta_{m-1}}+1,\omega^{\beta_m}] \text{ for } m\geq 2. 
\]

In either case, the intervals $(J_m)_{m\in\en}$ are pairwise disjoint,
compact and open in $[0,\omega^\alpha)$, and
\[
   [0,\omega^\alpha)=\bigcup_{m\in\en}J_m,
   \qquad
   \lim_{m\to\infty}\min J_m=\omega^\alpha.
\]
For every $\gamma\in\Gamma$ put
\[
   I_\gamma:=
   \bigcup_{\substack{m\in\en\\ \sigma_m=\gamma}}J_m.
\]
Clearly, the sets $\{I_\gamma\setsep\gamma\in\Gamma\}$ are pairwise
disjoint, their union is $[0,\omega^\alpha)$, and every $I_\gamma$ is open.

Fix $\gamma\in\Gamma$. Since $\gamma$ occurs in $(\sigma_m)$ infinitely
many times and $\min J_m\to\omega^\alpha$, the point $\omega^\alpha$ is
a limit point of $I_\gamma$. On the other hand, if
$\xi\in[0,\omega^\alpha)\setminus I_\gamma$, then $\xi\in J_m$ for some
$m$ with $\sigma_m\neq\gamma$, and the open set $J_m$ is disjoint from
$I_\gamma$. Consequently, $\overline{I_\gamma}=I_\gamma\cup\{\omega^\alpha\}$.

Now, let us identify the homeomorphism type of $I_\gamma$. Let
$m_1<m_2<\cdots$ be all the indices such that
$\sigma_{m_i}=\gamma$. If $\alpha=\beta+1$, then $I_\gamma$ is a disjoint
union of infinitely many intervals of the form $[\omega^\beta\cdot l+1,\omega^\beta\cdot(l+1)]$, possibly together with the initial interval $[0,\omega^\beta]$.
Therefore, the Mazurkiewicz--Sierpi\'nski classification gives that
$\overline{I_\gamma}$ is homeomorphic to $[0,\omega^{\beta+1}]$. If $\alpha$ is limit, then each $J_{m_i}$ is homeomorphic to
$[0,\omega^{\beta_{m_i}}]$ and $\lim_{i\to\infty}\beta_{m_i}=\alpha$. Again, the Mazurkiewicz--Sierpi\'nski classification gives that
$\overline{I_\gamma}$ is homeomorphic to $[0,\omega^\alpha]$.
In both cases, $\big(\overline{I_\gamma}\big)^{(\alpha)}
   =\{\omega^\alpha\}$, so the homeomorphism may be chosen to map $\omega^\alpha$ onto
$\omega^\alpha$. It follows that $I_\gamma$ is homeomorphic to
$[0,\omega^\alpha)$. This proves condition
\ref{it:splitIntoIntervals}.

Finally, fix $\eta<\omega^\alpha$. Choose $m_\eta\in\en$ such that $\min J_m>\eta$ whenever $m>m_\eta$. If $\min I_\gamma\leq\eta$, then $\sigma_m=\gamma$ for some
$m\leq m_\eta$. Hence $\{\gamma\in\Gamma\setsep\min I_\gamma\leq\eta\}
   \subset
   \{\sigma_1,\ldots,\sigma_{m_\eta}\}$,
and so this set is finite. This proves condition
\ref{it:separateIn}. In particular, every one-to-one sequence
$(\gamma_n)$ in $\Gamma$ eventually avoids this finite set, and therefore $\lim_{n\to\infty}\min I_{\gamma_n}=\omega^\alpha$.
\end{proof}

\begin{prop}
\label{prop:Construction with better bound}
Let $\alpha$ be a countable nonzero ordinal and $k \in \en$ with $k\geq 2$. Then
\[
d_{BM}^+\Big(\C(\omega^{\alpha}\cdot k),\C(\omega^{\alpha})\Big)\leq k + \sqrt{(k-1)(k+3)}.
\]
\end{prop}
\begin{proof}
Apply Lemma~\ref{lemma:split into subintervals} with
$\Gamma=\{1,\ldots,k\}$, and let $I_1,\ldots,I_k$ be the resulting
partition of $[0,\omega^\alpha)$. Put
\[
K:=[0,\omega^\alpha\cdot k],\qquad
L:=[0,\omega^\alpha+(k-1)],\qquad
e_j:=\omega^\alpha\cdot j,
\]
and, for $j=1,\ldots,k$, let
\[
J_1:=[0,\omega^\alpha],\qquad
J_j:=[\omega^\alpha\cdot(j-1)+1,\omega^\alpha\cdot j]
\quad (j\geq2).
\]
Thus $J_1,\ldots,J_k$ form a clopen partition of $K$, and each $J_j$ is
homeomorphic to $[0,\omega^\alpha]$, with $e_j$ corresponding to
$\omega^\alpha$.

For every
$i=1,\ldots,k$ pick a homeomorphism
$r_i:[0,\omega^\alpha]\longrightarrow \overline{I_i}$ such that $r_i(\omega^\alpha)=\omega^\alpha$.
For every $j=1,\ldots,k$, choose also a homeomorphism
$s_j:[0,\omega^\alpha]\to J_j$ such that
$s_j(\omega^\alpha)=e_j$, and put
\[
p_{i,j}:=s_j\circ r_i^{-1}:
\overline{I_i}\longrightarrow J_j.
\]
Then $p_{i,j}(\omega^\alpha)=e_j$ and, for every admissible choice of
indices,
\begin{equation}
\label{eq:common-coordinate-prop-better-bound}
p_{i,j}\circ p_{i,l}^{-1}
=s_j\circ s_l^{-1}
=p_{m,j}\circ p_{m,l}^{-1}.
\end{equation}

For $j=1,\ldots,k-1$, put $a_j:=\omega^\alpha+j\in L$. Fix
$\lambda=(\lambda_1,\ldots,\lambda_k)\in\Delta_k$. For $f\in\C(K)$
define $Tf\in\C(L)$ by
\begin{align*}
Tf(a_j)&:=f(e_j),
&&j=1,\ldots,k-1,\\
Tf(x)&:=\sum_{j=i}^k\lambda_j f(p_{i,j}(x))
       +\sum_{j=1}^{i-1}\lambda_j f(e_j),
&&x\in\overline{I_i},\quad i=1,\ldots,k.
\end{align*}
The formulas on the sets $\overline{I_i}$ agree at their common point
$\omega^\alpha$, where they all have the value
$\sum_{j=1}^k\lambda_jf(e_j)$. For every $i=1,\ldots,k$, the restriction of $Tf$ to
$\overline{I_i}$ is continuous, since it is a linear combination of
the continuous functions $f\circ p_{i,j}$ and constant functions.
Moreover, the sets $\overline{I_1},\ldots,\overline{I_k}$ are closed,
they cover $[0,\omega^\alpha]$, and for $i\neq l$ we have $\overline{I_i}\cap\overline{I_l}=\{\omega^\alpha\}$. At this common point all the formulas give the same value $\sum_{j=1}^k\lambda_j f(e_j)$. Consequently, $Tf$ is continuous on $[0,\omega^\alpha]$. Since the
points $a_1,\ldots,a_{k-1}$ are isolated, $Tf$ is continuous on $L$. It is clear that $T$ is linear and positive and that
$\norm{T}\leq1$.

We define $S:\C(L)\to\C(K)$ by
\begin{align*}
Sh(y)&:=h(a_l)+\frac{1}{\lambda_l}
 \big(
   h(p_{l,l}^{-1}(y))
   -h(p_{l+1,l}^{-1}(y))
 \big),
&&y\in J_l,\quad l=1,\ldots,k-1,\\
Sh(y)&:=\frac{1}{\lambda_k}
 \left(
   h(p_{k,k}^{-1}(y))
   -\sum_{j=1}^{k-1}\lambda_jh(a_j)
 \right),
&&y\in J_k.
\end{align*}
Since the sets $J_l$ are clopen and all the maps occurring in these
formulas are continuous, $Sh$ is continuous. Moreover,
\[
\norm{S}\leq
\max\left\{
   \frac{2}{\lambda_k}-1,\,
   1+\frac{2}{\lambda_l}\setsep l=1,\ldots,k-1
\right\}.
\]

Let us verify that $S=T^{-1}$. If $f\in\C(K)$, $l<k$, and $y\in J_l$,
then \eqref{eq:common-coordinate-prop-better-bound} gives
\[
Tf(p_{l,l}^{-1}(y))
-Tf(p_{l+1,l}^{-1}(y))
=
\lambda_l\big(f(y)-f(e_l)\big).
\]
For $y\in J_k$ we similarly have
\[
Tf(p_{k,k}^{-1}(y))
=
\lambda_kf(y)+\sum_{j=1}^{k-1}\lambda_jf(e_j).
\]
Together with $Tf(a_j)=f(e_j)$, these identities imply $STf=f$.

Conversely, let $h\in\C(L)$. We have
$TSh(a_j)=Sh(e_j)=h(a_j)$ for $j<k$. Fix
$i\in\{1,\ldots,k\}$ and $x\in\overline{I_i}$, and put
\[
   t:=r_i^{-1}(x),
   \qquad
   x_j:=r_j(t),
   \qquad j=1,\ldots,k.
\]
Then $x_i=x$ and, for every $j=1,\ldots,k$,
\[
   p_{i,j}(x)
   =
   s_j(r_i^{-1}(x))
   =
   s_j(t).
\]
In particular, for $j=1,\ldots,k-1$,
\begin{align*}
   p_{j,j}^{-1}(p_{i,j}(x))
   &=
   (r_j\circ s_j^{-1})(s_j(t))
   =
   r_j(t)
   =
   x_j,\\
   p_{j+1,j}^{-1}(p_{i,j}(x))
   &=
   (r_{j+1}\circ s_j^{-1})(s_j(t))
   =
   r_{j+1}(t)
   =
   x_{j+1}.
\end{align*}
Hence, for $j=i,\ldots,k-1$, the definition of $S$ gives
\begin{align*}
   Sh(p_{i,j}(x))
   &=
   h(a_j)
   +
   \frac{1}{\lambda_j}
   \Big(
      h(p_{j,j}^{-1}(p_{i,j}(x)))
      -
      h(p_{j+1,j}^{-1}(p_{i,j}(x)))
   \Big)\\
   &=
   h(a_j)
   +
   \frac{h(x_j)-h(x_{j+1})}{\lambda_j}.
\end{align*}
Similarly,
\[
   p_{k,k}^{-1}(p_{i,k}(x))
   =
   (r_k\circ s_k^{-1})(s_k(t))
   =
   r_k(t)
   =
   x_k,
\]
and so
\[
   Sh(p_{i,k}(x))
   =
   \frac{1}{\lambda_k}
   \left(
      h(x_k)
      -
      \sum_{j=1}^{k-1}\lambda_jh(a_j)
   \right).
\]

Using these formulas and $Sh(e_j)=h(a_j)$, we obtain
\begin{align*}
TSh(x)
&=
\sum_{j=i}^{k-1}
 \big(
   \lambda_jh(a_j)+h(x_j)-h(x_{j+1})
 \big)
 +h(x_k)-\sum_{j=1}^{k-1}\lambda_jh(a_j)
 +\sum_{j=1}^{i-1}\lambda_jh(a_j)\\ 
 &=h(x_i)=h(x).
\end{align*}
Hence $TS=\Id$ and $S=T^{-1}$.

We have therefore shown that, for every $\lambda\in\Delta_k$,
\[
\norm{T}\cdot\norm{T^{-1}}
\leq
\max\left\{
   \frac{2}{\lambda_k}-1,\,
   1+\frac{2}{\lambda_l}\setsep l=1,\ldots,k-1
\right\}.
\]
Since $\Delta_k$ is invariant under permutations of the coordinates,
Lemma~\ref{lemma:the constant} gives a choice of $\lambda$ for which the
right-hand side equals
\[
C(k)=k+\sqrt{(k-1)(k+3)}.
\]
Finally, $L$ is homeomorphic to $[0,\omega^\alpha]$ by the
Mazurkiewicz--Sierpi\'nski classification, which proves the assertion.
\end{proof}

\begin{lemma}
\label{lemma:orindal_repre}
Let $\alpha, \beta$ be nonzero countable ordinals and $n \in \en$.
\begin{enumerate}[label=(\roman*)]
    \item\label{it:ordinalFormWithBeta} Each ordinal $y \in [1, \omega^{\alpha\cdot n}\cdot\beta]$ can be uniquely expressed in the form 
$y=\sum_{i=1}^k \omega^{\alpha\cdot(n-i+1)}\cdot\gamma_i$,
for some $1 \leq k \leq n+1$, $\gamma_1 \in [0, \beta]$, $\gamma_i \in [0, \omega^{\alpha})$ for $2 \leq i \leq k$, and such that $\gamma_k \neq 0$.  
    \item\label{it:ordinalFormWithoutBeta} 
Each ordinal $y \in [1, \omega^{\alpha\cdot n})$ can be uniquely expressed in the form 
$y=\sum_{i=1}^k \omega^{\alpha\cdot(n-i)}\cdot\gamma_i$,
for some $1 \leq k \leq n$, $\gamma_i \in [0, \omega^{\alpha})$ for $1 \leq i \leq k$, and such that $\gamma_k \neq 0$.
\end{enumerate}
\end{lemma}

\begin{proof}[Sketch of the proof]
The proof is easy and standard (even though with full details quite tedious), so we give just a sketch of the proof and leave the details to the interested reader.

\smallskip

\noindent\ref{it:ordinalFormWithBeta}: We write $y$ in its Cantor normal form and decompose further as
\[
y = \underbrace{\omega^{\alpha_1}\cdot n_1 + \ldots + \omega^{\alpha_{{k_1}}}\cdot n_{k_1}}_{=:y_1} + \ldots + \underbrace{\omega^{\alpha_{k_n+1}}\cdot n_{k_n+1} + \ldots + \omega^{\alpha_{{k_{n+1}}}}\cdot n_{k_{n+1}}}_{=:y_{n+1}},
\]
where
\[
\alpha_1>\ldots>\alpha_{k_1}\geq \alpha\cdot n > \alpha_{k_1+1}>\ldots>\alpha_{k_i}\geq \alpha\cdot(n-i+1)>\alpha_{k_i+1}>\ldots>\alpha_{k_{n+1}}\geq 0
\]
and $n_i\in\en\cup\{0\}$ for $i=1,\ldots,k_{n+1}$. Next, we find the unique $\gamma_1\in [0,\beta]$ satisfying $y_1 = \omega^{\alpha\cdot n}\cdot \gamma_1$ and for every $2\leq i\leq n+1$ we find the unique $\gamma_i\in[0,\omega^{\alpha})$ with $y_i=\omega^{\alpha\cdot (n-i+1)}\cdot \gamma_i$. Finally, we pick $k$ such that $y_k\neq 0$ and $y_i=0$ for $i > k$.

\smallskip

\noindent\ref{it:ordinalFormWithoutBeta}: Is similar to the proof of \ref{it:ordinalFormWithBeta} and therefore left to the reader.
\end{proof}

\begin{prop}
\label{prop:constructing-isomorphism-new-attempt}
Let $\alpha, \beta$ be countable nonzero ordinals and $n \in \en$. Then
\begin{itemize}
    \item[(i)] $
d_{BM}^+\Big(\C( \omega^{\alpha}+\beta),\C( \omega^{\alpha\cdot n}\cdot\beta)\Big)\leq n+1+\sqrt{n(n+4)},$
\item[(ii)] $
d_{BM}^+\Big(\C(\omega^{\alpha}),\C( \omega^{\alpha\cdot n})\Big)\leq n+\sqrt{(n-1)(n+3)}
.$
\end{itemize}
\end{prop}

\begin{proof}We prove the statement (i) in detail, while we only indicate the proof of (ii),
because it follows by just a minor modification of the proof of (i).

To start with, we consider the set 
\[
\begin{split}
\Theta=\bigl\{(\gamma_1,\ldots,\gamma_k):1\leq k\leq n+1,\;
&\gamma_1<\beta
\text{ and }
\gamma_i<\omega^\alpha\text{ for }2\leq i\leq k,\\
&\text{or }
\gamma_1=\beta
\text{ and }
\gamma_2=\cdots=\gamma_k=0
\bigr\}.
\end{split}
\]
and given $\gamma=(\gamma_1,\ldots,\gamma_k)\in\Theta$ we define $l(\gamma):=k$ to be the length of the tuple $\gamma$, if $k=1$ we write rather $\gamma_1\in\Theta$ instead of $(\gamma_1)\in\Theta$. Given $\gamma\in\Theta$ and $1\leq i<l(\gamma)$ we put $\gamma|_i:=(\gamma_1,\ldots,\gamma_i)$. Further, given $\gamma\in\Theta$ we let $y_\gamma:=\sum_{i=1}^{l(\gamma)}\omega^{\alpha\cdot(n-i+1)}\cdot\gamma_i$. Put $\Theta':=\{\gamma\in\Theta\setsep \gamma_{l(\gamma)}\neq 0\}$ and recall that by Lemma~\ref{lemma:orindal_repre} the set $\Theta'$ is bijectively mapped onto $[1,\omega^{\alpha\cdot n}\cdot\beta]$ using the map $\gamma\mapsto y_\gamma$. On the set $\Theta'$ we shall consider the ordering and topology inherited from $[1,\omega^{\alpha\cdot n}\cdot\beta]$ which means that $\Theta'$ is ordered lexicographically and endowed with the order topology. Let us shortly mention what basic open neighborhoods of $\gamma = (\gamma_1,\ldots,\gamma_k)\in\Theta'$ look like. If $\gamma_k$ is a limit ordinal, basic clopen neighborhoods of $\gamma$ are of the form
\[
U_\delta(\gamma):=\{\gamma\}\cup\{(\gamma_1,\ldots,\gamma_{k-1},\eta_1,\ldots,\eta_l)\in\Theta'\setsep \delta<\eta_1<\gamma_k\},
\]
where $\delta\in [0,\gamma_k)$. If $\gamma_k$ is a successor ordinal and $k=n+1$ then $\gamma$ is an isolated point and if $\gamma_k$ is a successor ordinal and $k\leq n$, then basic clopen neighborhoods of $\gamma$ are of the form
\[
U_\delta(\gamma):=\{\gamma\}\cup\{(\gamma_1,\ldots,\gamma_{k-1},\gamma_k - 1,\eta_1,\ldots,\eta_l)\in\Theta'\setsep \eta_1 > \delta\},
\]
where $\delta\in [0,\omega^\alpha)$.

Note that if we add one isolated point to $\Theta'$, denoted as $*$, then $\Theta'\cup\{*\}$ is still homeomorphic to $[1,\omega^{\alpha\cdot n}\cdot\beta]$. Thus, it suffices to construct positive isomorphism $T:\C([0,\omega^\alpha + \beta])\to \C(\Theta'\cup\{*\})$ with $\norm{T}\leq 1$ and $\norm{T^{-1}}\leq C(n+1)$ (see Lemma~\ref{lemma:the constant}). In order to shorten notation we shall denote $X:=[0,\omega^\alpha + \beta]$ and $N:=\Theta'\cup\{*\}$.

Let us first define a bijection $\rho:N \rightarrow X$ as follows. Let 
\[
\mathcal P:=\{\gamma\in\Theta:
l(\gamma)\leq n,\ \gamma_1<\beta\}.\]
We regard $\mathcal P$ merely as an index set. Apply
Lemma~\ref{lemma:split into subintervals} with
$\Gamma=\mathcal P$, and denote the resulting partition of
$[0,\omega^\alpha)$ by
$\{I_\eta:\eta\in\mathcal P\}$. For every $\eta\in\mathcal P$, choose a homeomorphism $\rho_\eta:[1,\omega^\alpha]\longrightarrow\overline{I_\eta}$ such that $\rho_\eta(\omega^\alpha)=\omega^\alpha$. Finally, we define $\rho:N \rightarrow X$ as follows. We put $\rho(*):=\omega^\alpha$ and 
\[
    \rho(\gamma)=
    \begin{cases}
        \omega^\alpha+\gamma_1,
            &l(\gamma)=1,\\[1mm]
        \rho_{\gamma|_{l(\gamma)-1}}
        \bigl(\gamma_{l(\gamma)}\bigr),
            &l(\gamma)>1.
    \end{cases}
\]
Note that then restriction of $\rho$ to each of the sets $\{(\gamma_1,\ldots,\gamma_{k-1},\gamma')\setsep \gamma'\in [1,\omega^\alpha)\}$ (with $k\geq 2$) is a homeomorphism onto $I_{(\gamma_1,\ldots,\gamma_{k-1})}$ and it easily follows that $\rho$ is indeed a bijection. Also, we can observe that restriction of $\rho$ to each of the sets $L_k:=\{\gamma\in\Theta'\setsep l(\gamma)=k\}$, $1\leq k\leq n+1$ is a homeomorphism. Moreover, $\rho(L_1)=[\omega^\alpha+1,\omega^\alpha+\beta]$ and for  $k\geq 2$ we have $\rho(L_k) =
\bigcup\{I_\eta\setsep \eta\in\mathcal P, l(\eta)=k-1\}$, so for any $k$ the set $\rho(L_k)$ is open in $X\setminus\{\omega^\alpha\}$.

For $\gamma\in\Theta'$, denoting for $1\leq i < l(\gamma)$ the sequence $(\gamma_1,\ldots,\gamma_i+1)$ as $\gamma|_i+1$, we define the parent of $\gamma$ by
\[
    \partial\gamma:=
    \begin{cases}
        *,&l(\gamma)=1,\\
        \gamma|_{l(\gamma)-1}+1,&l(\gamma)>1.
    \end{cases}
\]
Thus, putting $L_0=\{*\}$, we have
$\partial(L_k)\subset L_{k-1}$.
It follows directly from the description of the topology on the
sets \(L_k\) that every restriction
\[
    \partial|_{L_k}:L_k\longrightarrow L_{k-1}
\]
is continuous (indeed, for $k=1$ this is a constant map and for $k\geq 2$ the mapping $\partial|_{L_k}$ is constant on each open set $\{\eta^\frown\xi:\xi\in[1,\omega^\alpha)\}$ for $\eta\in \mathcal P$ with $l(\eta)=k-1$).

For $i=1,\ldots,n+1$, define
$R_i:N\to X$ by $R_i(*)=\omega^\alpha$ and, for $\gamma\in L_k$, by
\[
    R_i(\gamma):=
    \begin{cases}
        \rho\bigl(\partial^{\,k-i}\gamma\bigr),&i\leq k,\\
        \omega^\alpha,&i>k,
    \end{cases}
\]
where as usual, $\partial^0:=\Id$. Equivalently,
\[
    R_i(\gamma)=
    \begin{cases}
        \rho(\gamma|_i+1),&i<l(\gamma),\\
        \rho(\gamma),&i=l(\gamma),\\
        \omega^\alpha,&i>l(\gamma).
    \end{cases}
\]

We claim that all the maps $R_i$ are continuous. Since $N$ and
$X$ are metrizable, it suffices to prove sequential continuity.
Moreover, given a sequence $(\gamma_j)$ converging to $\gamma$,
it is enough to show that every subsequence of
$(R_i(\gamma_j))$ has a further subsequence converging to
$R_i(\gamma)$. Thus, let $\gamma_j\to\gamma$. The case $\gamma=*$ is immediate,
since $*$ is isolated. Given any subsequence of $(\gamma_j)$, we
may pass to a further subsequence such that $l(\gamma_j)=k$ is
constant. The description of the lexicographic topology gives $k\geq l(\gamma)$. If $k=l(\gamma)$, then
$R_i(\gamma_j)\longrightarrow R_i(\gamma)$ follows from the continuity of $R_i|_{L_k}$. We may therefore
assume that $k>l(\gamma)$. In this case, the description of the
basic neighbourhoods of $\gamma$ gives
$\gamma_j|_{l(\gamma)}+1\longrightarrow\gamma$. It follows that $R_i(\gamma_j)\longrightarrow R_i(\gamma)$ for $i\leq l(\gamma)$. If $i>k$, then $R_i(\gamma_j)=\omega^\alpha=R_i(\gamma)$. If $l(\gamma) < i\leq k$, put $\delta_j:=\partial^{\,k-i}\gamma_j\in L_i$. Then $R_i(\gamma_j)=\rho(\delta_j)$. After passing to a further
subsequence, the points $R_i(\gamma_j)$ either belong
to pairwise distinct sets $I_\eta$, in which case $ R_i(\gamma_j)\longrightarrow\omega^\alpha = R_i(\gamma)$ by Lemma~\ref{lemma:split into subintervals}\ref{it:separateIn}, or they
belong to a fixed $I_\eta$. In the latter case, we may write $\delta_j=\eta^\frown\xi_j$ and hence $R_i(\gamma_j)=\rho_\eta(\xi_j)$. Since $\gamma_j\to\gamma$, the description of
the basic neighbourhoods of $\gamma$ implies
$\xi_j\to\omega^\alpha$ and therefore $R_i(\gamma_j)
    =\rho_\eta(\xi_j)
    \longrightarrow
    \rho_\eta(\omega^\alpha)
    =\omega^\alpha = R_i(\gamma)$. We have shown that every subsequence has a further subsequence with
the required convergence. Hence $R_i(\gamma_j)\to R_i(\gamma)$,
and $R_i$ is continuous.

We shall also need the following topological consequence of the
construction: whenever sequence $(\gamma_j)$ in $N\setminus\{*\}$ satisfies
$\rho(\gamma_j)\to \omega^\alpha$ and
$(\gamma_{j_k},\partial \gamma_{j_k})
    \to (\gamma,\eta)\in N\times N$ 
for some subsequence $(\gamma_{j_k})$, then $\gamma=\eta$.

Indeed, after replacing the sequence by the indicated subsequence
and then passing to a further subsequence, we may assume that
$l(\gamma_j)=k$ is constant and that each coordinate sequence is
either constant or strictly increasing. Necessarily $k\geq 2$,
since
\[
    \rho(L_1)=[\omega^\alpha+1,\omega^\alpha+\beta]
\]
does not accumulate at $\omega^\alpha$. If the first $k-1$ coordinates are constant, then
$\rho(\gamma_j)\to\omega^\alpha$ implies
$\gamma_j(k)\to\omega^\alpha$. Hence $\gamma_j$ and the constant
sequence $(\partial\gamma_j)$ converge to the same point, so
$\gamma=\eta$. Otherwise, after passing to a further subsequence, the lexicographic
ordering gives
\[
    \gamma_j<\partial\gamma_j<\gamma_{j+1}.
\]
It follows that \(\gamma_j\) and \(\partial\gamma_j\) again have the
same limit. Thus \(\gamma=\eta\), proving the claim.

We now isolate the operator-theoretic part of the argument, which
will also be used in the proof of (ii).

\begin{claim}[Finite-level coding principle]
Let \(X\) and \(N\) be compact Hausdorff spaces, let \(p\in X\) and
\(*\in N\), and suppose that
\[
    N=\{*\}\sqcup L_1\sqcup\cdots\sqcup L_m.
\]
Put \(L_0=\{*\}\), and let $\partial:\bigcup_{k=1}^mL_k\longrightarrow N$ be such that $\partial(L_k)\subset L_{k-1}$ and \(\partial|_{L_k}\) is continuous for every \(k\). Suppose that there is a bijection $\rho:N\longrightarrow X$ with $\rho(*)=p$ such that every \(\rho(L_k)\) is open in \(X\setminus\{p\}\) and $\rho|_{L_k}:L_k\longrightarrow\rho(L_k)$ is a homeomorphism.

For $i=1,\ldots,m$, define $R_i:N\to X$ by $R_i(*)=p$ and, for $x\in L_k$ with $k=1,\ldots,m$, by
\[
    R_i(x)=
    \begin{cases}
        \rho\bigl(\partial^{\,k-i}x\bigr),&i\leq k,\\
        p,&i>k.
    \end{cases}
\]
Assume that all the maps \(R_i\) are continuous and that the following
condition holds: whenever sequence $(x_j)$ in $N\setminus\{*\}$ satisfies
$\rho(x_j)\to p$ and
$(x_{j_k},\partial x_{j_k})
    \to (x,y)\in N\times N$ 
for some subsequence $(x_{j_k})$, then $x=y$.

Then \[d_{BM}^+(\C(X),\C(N))\leq C(m).\]
\end{claim}

\begin{proof}[Proof of the claim]
Fix $\lambda=(\lambda_1,\ldots,\lambda_m)\in\Delta_m$ and define
\[Tf(x):=\sum_{i=1}^m\lambda_i f(R_i(x)),\qquad x\in N.\]
Since the
maps $R_i$ are continuous, $T$ is a well-defined positive norm-one
operator.

Define $S:\C(N)\to\C(X)$ by $Sh(p):=h(*)$ and
\begin{equation}
\label{eq:finite-level-inverse}
    Sh(\rho(x))
    :=
    h(*)+\frac{h(x)-h(\partial x)}{\lambda_k},
    \qquad x\in L_k,\quad k=1,\ldots,m.
\end{equation}
Since $\rho$ is a bijection, this defines $Sh$ on all of $X$. The
assumptions on $\rho|_{L_k}$ and $\partial|_{L_k}$ show that $Sh$ is
continuous on every open set $\rho(L_k)$, and hence on
$X\setminus\{p\}$.

We verify continuity at $p$. Let $(x_j)\subset N\setminus\{*\}$
satisfy $\rho(x_j)\to p$. We claim that
$h(x_j)-h(\partial x_j)\to0$. Otherwise, after passing to a
subsequence, there would be $\varepsilon>0$ such that
$|h(x_j)-h(\partial x_j)|\geq\varepsilon$ for every $j\in\en$. By
compactness of $N\times N$, we may pass to a further subsequence such
that $(x_j,\partial x_j)\to(x,y)$ for some $x,y\in N$. By assumption,
$x=y$, which contradicts the continuity of $h$. Therefore
$h(x_j)-h(\partial x_j)\to0$. Since the set
$\{\lambda_1,\ldots,\lambda_m\}$ is finite,
\eqref{eq:finite-level-inverse} gives
$Sh(\rho(x_j))\to h(*)=Sh(p)$. Thus $Sh\in\C(X)$.

If $x\in L_1$, then $\partial x=*$, and hence
\[
    Sh(\rho(x))
    =
    \frac{1}{\lambda_1}h(x)
    +
    \left(1-\frac{1}{\lambda_1}\right)h(*).
\]
For $x\in L_k$ with $k\geq2$, we use
\eqref{eq:finite-level-inverse} directly. Consequently,
\begin{equation}
\label{eq:finite-level-inverse-bound}
    \norm{S}
    \leq
    \max\left\{
        \frac{2}{\lambda_1}-1,\,
        \max_{2\leq k\leq m}
        \left(1+\frac{2}{\lambda_k}\right)
    \right\},
\end{equation}
where the second maximum is omitted when $m=1$.

For $x\in L_k$, the definition of the maps $R_i$ gives
$R_i(x)=R_i(\partial x)$ for $i<k$ and
$R_i(\partial x)=p$ for $i\geq k$. It follows that
\begin{equation}
\label{eq:finite-level-parent-identity}
    Tf(x)-Tf(\partial x)
    =
    \lambda_k\bigl(f(\rho(x))-f(p)\bigr).
\end{equation}
Since $Tf(*)=f(p)$, we obtain
\[
    S(Tf)(\rho(x))
    =
    Tf(*)+
    \frac{Tf(x)-Tf(\partial x)}{\lambda_k}=
    f(p)+f(\rho(x))-f(p)
    =
    f(\rho(x)).
\]
Moreover, $S(Tf)(p)=Tf(*)=f(p)$, and hence
$ST=\Id_{\C(X)}$.

Conversely, $T(Sh)(*)=Sh(p)=h(*)$. Suppose inductively that
$T(Sh)(\partial x)=h(\partial x)$ for some $x\in L_k$. Applying
\eqref{eq:finite-level-parent-identity} to $Sh$, we obtain
\[
    T(Sh)(x)
    =
    T(Sh)(\partial x)
    +\lambda_k\bigl(Sh(\rho(x))-Sh(p)\bigr)=
    h(\partial x)+h(x)-h(\partial x)
    =
    h(x).
\]
Induction over the levels gives $TS=\Id_{\C(N)}$. Thus $S=T^{-1}$.
By \eqref{eq:finite-level-inverse-bound} and
Lemma~\ref{lemma:the constant}, we may choose
$\lambda\in\Delta_m$ so that
$\norm{T}\norm{T^{-1}}\leq C(m)$. This proves the claim.
\end{proof}

We apply the claim to the model constructed above with
\[
    p=\omega^\alpha
    \qquad\text{and}\qquad
    m=n+1.
\]
All its assumptions have been verified, therefore
\[
d_{BM}^+\Bigl(
    \C(\omega^\alpha+\beta),
    \C(\omega^{\alpha\cdot n}\cdot\beta)
\Bigr)
\leq
C(n+1)
=
n+1+\sqrt{n(n+4)}.
\]
This finishes the proof of (i).

In order to prove (ii), we adapt the preceding construction as
follows. First, consider the set of tuples representing the ordinal
decomposition:
\[
    \Theta
    :=
    \{\omega^{\alpha\cdot n}\}
    \cup
    \bigcup_{k=1}^n[0,\omega^\alpha)^k.
\]
For every tuple
$\gamma=(\gamma_1,\ldots,\gamma_k)\in\Theta
\setminus\{\omega^{\alpha\cdot n}\}$, let $l(\gamma):=k$ denote its
length. We denote the truncation of $\gamma$ to its first $i$
coordinates by
$\gamma|_i:=(\gamma_1,\ldots,\gamma_i)$ and adopt the convention
$\gamma|_0:=\varnothing$. Put $\Theta'
    :=
    \{\omega^{\alpha\cdot n}\}
    \cup
    \left\{
        \gamma\in\Theta\setminus\{\omega^{\alpha\cdot n}\}:
        \gamma_{l(\gamma)}\neq0
    \right\}$. For $\gamma\neq\omega^{\alpha\cdot n}$, let $y_\gamma
    :=
    \sum_{i=1}^{l(\gamma)}
    \omega^{\alpha\cdot(n-i)}\cdot\gamma_i$, and put
$y_{\omega^{\alpha\cdot n}}:=\omega^{\alpha\cdot n}$. By
Lemma~\ref{lemma:orindal_repre}\ref{it:ordinalFormWithoutBeta}, the
map $\gamma\mapsto y_\gamma$ identifies $\Theta'$ with
$[1,\omega^{\alpha\cdot n}]$. We equip $\Theta'$ with the
corresponding order topology. In particular, $\Theta'$ is
homeomorphic to $[0,\omega^{\alpha\cdot n}]$.

Let
\[
    \mathcal P_0
    :=
    \{\varnothing\}
    \cup
    \bigcup_{k=1}^{n-1}[0,\omega^\alpha)^k.
\]
Apply Lemma~\ref{lemma:split into subintervals} with
$\Gamma=\mathcal P_0$, and denote the resulting partition of
$[0,\omega^\alpha)$ by
$\{I_\eta:\eta\in\mathcal P_0\}$. For every
$\eta\in\mathcal P_0$, fix a homeomorphism $\rho_\eta:[1,\omega^\alpha]
    \longrightarrow\overline{I_\eta}$
such that $\rho_\eta(\omega^\alpha)=\omega^\alpha$. Define a bijection $\rho:\Theta'\to[0,\omega^\alpha]$ by
$\rho(\omega^{\alpha\cdot n}):=\omega^\alpha$ and, for $\gamma\neq\omega^{\alpha\cdot n}$, by
\[
    \rho(\gamma)
    :=
    \rho_{\gamma|_{l(\gamma)-1}}
    \bigl(\gamma_{l(\gamma)}\bigr).
\]
Put
\[
    L_0:=\{\omega^{\alpha\cdot n}\},
    \qquad
    L_k:=\{\gamma\in\Theta':l(\gamma)=k\},
    \quad k=1,\ldots,n,
\]
and define
\[
    \partial\gamma:=
    \begin{cases}
        \omega^{\alpha\cdot n},&l(\gamma)=1,\\
        \gamma|_{l(\gamma)-1}+1,&l(\gamma)>1.
    \end{cases}
\]
Similar arguments as in (i) show that assumptions of the Claim above are verified with
\[
    X=[0,\omega^\alpha],\qquad
    N=\Theta',\qquad
    p=\omega^\alpha,\qquad
    m=n,\qquad *=\omega^{\alpha\cdot n}.
\]
We obtain
\[
d_{BM}^+\Bigl(
    \C(\omega^\alpha),
    \C(\omega^{\alpha\cdot n})
\Bigr)
\leq
C(n)
=
n+\sqrt{(n-1)(n+3)}.
\]
This proves (ii).
\end{proof}

Next, before we proceed to the proof of Theorem \ref{thm:class}, we need the following lemma.

\begin{lemma}
\label{lemma:another form of ordinal, version 2}
Let $\alpha, \beta$ be countable nonzero ordinals such that $\alpha \leq \beta<\alpha\cdot\omega$, and let $m \in \en$. Then, there exists $n \in \en$ and a countable ordinal $\gamma<\omega^{\alpha}$ such that $[1, \omega^\beta\cdot m]$ is homeomorphic to $[1, \omega^{\alpha\cdot n}\cdot\gamma]$.
\end{lemma}

\begin{proof}

We find a maximal $n \in \en$ such that $\alpha\cdot n \leq \beta$. Further, we find an ordinal $\delta$ such that $\alpha\cdot n+\delta=\beta$ (yes, this is possible, see e.g. \cite[Lemma 2.25 (ii)]{Jechsbook}), and we set $\gamma=\omega^{\delta}m$. Then, 
\[\omega^{\alpha\cdot n}\cdot\gamma=\omega^{\alpha\cdot n}\cdot(\omega^{\delta}\cdot m)=\omega^{\alpha\cdot n +\delta}\cdot m=\omega^{\beta}\cdot m.\]

Now, assuming that $\gamma \geq \omega^{\alpha}$, we would have $\omega^{\beta}\cdot m=\omega^{\alpha\cdot n}\cdot\gamma \geq \omega^{\alpha\cdot(n+1)}$, 
which implies $\beta \geq \alpha\cdot(n+1)$, and this contradicts the maximality of $n$. 
\end{proof}

Now, we are ready to prove Theorem \ref{thm:class}.

\begin{proof}[Proof of Theorem \ref{thm:class}.]
The implication \ref{it:positivelyIso} $\Rightarrow$ \ref{it:isoAndHeight} follows directly from Theorem~\ref{thm:embedMain}. The equivalence between \ref{it:isoAndHeight} and \ref{it:height} is a consequence of the classical classification of $\C(K)$ spaces by Bessaga and Pe{\l}czy\'nski (see Remark \ref{rem:szlenkAndIso}).

We prove the implication \ref{it:height} $\Rightarrow$ \ref{it:positivelyIso}. Let $K$ and $L$ be infinite countable compact spaces such that $\height(K) \leq \height(L) < \Gamma(\height(K))$. By the Mazurkiewicz-Sierpi\'nski classification, we may assume without loss of generality that
\[
K = [1, \omega^\alpha \cdot k] \quad \text{and} \quad L = [1, \omega^\beta \cdot m],
\]
where $\alpha, \beta$ are countable ordinals and $k, m \in \en$.
The condition on the heights implies that $\alpha \leq \beta < \alpha\cdot\omega$.

We construct a positive isomorphism between $\C(K)$ and $\C(L)$ via a chain of positive isomorphisms.
First, by Proposition \ref{prop:Construction with better bound}, there exists a positive isomorphism
\[
\C([1, \omega^\alpha \cdot k]) \simeq \C([1, \omega^\alpha]).
\]
Next, since $\alpha \leq \beta < \alpha\cdot\omega$, we apply Lemma \ref{lemma:another form of ordinal, version 2} to find $n \in \en$ and $\gamma < \omega^\alpha$ such that $[1, \omega^\beta \cdot m]$ is homeomorphic to $[1, \omega^{\alpha\cdot n}\cdot\gamma]$. This homeomorphism induces a positive isomorphism (an isometry)
\[
\C([1, \omega^\beta \cdot m]) \simeq \C([1, \omega^{\alpha\cdot n}\cdot\gamma]).
\]
It remains to find a positive isomorphism between $\C([1, \omega^\alpha])$ and $\C([1, \omega^{\alpha\cdot n}\cdot\gamma])$.
Since $\gamma < \omega^\alpha$, the space $[1, \omega^\alpha]$ is homeomorphic to $[1, \omega^\alpha + \gamma]$. Thus, it suffices to construct a positive isomorphism
\[
T: \C([1, \omega^\alpha + \gamma]) \to \C([1, \omega^{\alpha\cdot n}\cdot\gamma]).
\]
The existence of such an operator is guaranteed by Proposition \ref{prop:constructing-isomorphism-new-attempt}(i).
Composing these operators yields the required positive isomorphism from $\C(K)$ to $\C(L)$.
\end{proof}

\begin{remark}\label{rem:upperBound}
The upper bound in Proposition \ref{prop:constructing-isomorphism-new-attempt}(ii) for the specific case $\alpha = 1$ can be derived from the proof of \cite[Theorem 1.4]{CandidoGalegoComega}. Although that result is formulated for the general Banach--Mazur distance, an inspection of the proof reveals that the constructed isomorphism is indeed positive. However, even for $\alpha=1$, the model considered in the proof of Proposition~\ref{prop:constructing-isomorphism-new-attempt} is more flexible: it relies on $n$ variables ($\lambda_i$) allowing for finer optimization, whereas the model in \cite{CandidoGalegoComega} depends on a single variable $A$. Moreover, the use of our tuple notation allows for a more compact presentation of both the isomorphism and its inverse. Finally, we emphasize that extending this construction to the case $\alpha>1$ required overcoming non-trivial topological obstacles not present in the $\alpha=1$ setting.
\end{remark}

Furthermore, by combining our construction with the main result of \cite{rondos-somaglia}, we obtain the following tight bounds on the distances between spaces of continuous functions on ordinal intervals. Observe that the difference between the upper and lower bounds is at most $2$ in all cases.

\begin{cor}
\label{cor:estimatesForOmegaToAlpha}
For each countable nonzero ordinal $\alpha$ and $n \in \en$,
\[
2n-1 \leq d_{BM}(\C(\omega^{\alpha}), \C(\omega^{\alpha\cdot n})) \leq d_{BM}^+(\C(\omega^{\alpha}), \C(\omega^{\alpha\cdot n})) \leq n+\sqrt{(n-1)(n+3)}.
\]
\end{cor}

\begin{proof}
The upper estimate is a direct consequence of Proposition \ref{prop:constructing-isomorphism-new-attempt}(ii). To deduce the lower estimate, consider the Cantor normal form of $\alpha$. Let $\omega^{\beta}\cdot k$ be the leading term of this expansion (where $k \in \en$ and $\beta$ is an ordinal). Then we have
\[
\height([1, \omega^{\alpha}]) = \alpha+1 \leq \omega^{\beta}\cdot(k+1).
\]
On the other hand,
\[
\height([1, \omega^{\alpha\cdot n}]) = \alpha\cdot n+1 > \omega^{\beta}\cdot nk.
\]
Consequently, by applying \cite[Theorem 1.1]{rondos-somaglia}, we obtain
\[
d_{BM}(\C(\omega^{\alpha}), \C( \omega^{\alpha\cdot n})) \geq \frac{2nk-k}{k}=2n-1.
\]
\end{proof}

\section{Quantitative estimates and upper constructions}\label{sec:estimates}

In this section, we improve estimates for the positive Banach--Mazur distance and collect upper constructions needed later. Our main result establishes the exact value of the distance for a specific class of ordinals. This result is new even for the classical Banach--Mazur distance. The upper bound is a direct consequence of Corollary~\ref{cor:estimatesForOmegaToAlpha}, so the primary contribution is the lower bound.

\begin{thm}
\label{thm:the exact values of the distances}
For each countable ordinal $\alpha$ and $n \in \en$,
\[
\begin{aligned}
d_{BM}(\C(\omega^{\omega^{\alpha}}), \C(\omega^{\omega^{\alpha}\cdot n}))
=d_{BM}^+(\C(\omega^{\omega^{\alpha}}), \C(\omega^{\omega^{\alpha}\cdot n}))
=n+\sqrt{(n-1)(n+3)}. 
\end{aligned}
\]
\end{thm}

\begin{remark}\label{rem:piasAnalogy}
For the case $\alpha = 0$, the lower bound in Theorem~\ref{thm:the exact values of the distances} constitutes the main result of the recent paper by Malec and Piasecki \cite{MaPi25}. Our proof of the lower bound adapts their argument, utilizing techniques from \cite{rondos-somaglia} to extend the result to this general ordinal setting.
\end{remark}

We also record a general lower estimate for compacta with finite last
derivatives and a two-branch upper construction valid at every countable rank.
The exact forward ordinal formula is proved later in
Theorem~\ref{thm:pd-forward-all-ranks}.

\begin{prop}\label{prop:pd-higher-two-upper}
For every countable nonzero ordinal $\alpha$,
\[
d_{BM}^+\bigl(\C(\omega^\alpha),\C(\omega^\alpha\cdot2)\bigr)
\leq2+\sqrt3.
\]
\end{prop}

In the remainder of this section we head towards the proofs of the results announced above.

\subsection{Proof of Theorem~\ref{thm:the exact values of the distances}}
We begin by proving Theorem~\ref{thm:lower estimate general isomorphism}. Combined with the results from Section~\ref{sec:main}, Theorem~\ref{thm:the exact values of the distances} will follow as a direct consequence. We first collect several auxiliary results.

\begin{lemma}  
\label{lemma:the lower estimate C(n)}
Let $K_1, K_2$ be compact spaces, $T: \C(K_1) \rightarrow \C(K_2)$ be an isomorphism with $\norm{T^{-1}}=1$, let $n \in \en, n \geq 2$, and let $\xi \in \er$. Assume that for each $\ep>0$, there exist norm-$1$ functions $f_1, \ldots, f_n, g, h \in \C(K_1, [0, 1])$ and a clopen set $U \subseteq K_1$ such that:
\begin{enumerate}[label=(a-\alph*)]
    \item\label{it:supports} for each $i=1, \ldots, n-1$, $\spt f_i \subseteq \{x \in K_1: f_{i+1}(x)=1\}$,
    \item\label{it:supports2} $\spt f_n \subseteq U \subseteq \{x \in K_1: g(x)=1\}$,
    \item\label{it:supports3} $U \cap \spt h=\emptyset$,
    \item\label{it:mapped supports} the sets $S_i=\{y \in K_2: \abs{Tf_i(y)} \geq \ep\}_{i=1}^n$ are pairwise disjoint,
    \item\label{it:almost 0} for each $i=1, \ldots, n-1$, $S_i \subseteq \{y \in K_2:\abs{T(T^{-1}(\mathbf{1}_{K_2})\cdot \chi_U)(y)} < \ep\}$,
    \item\label{it:almost constant} for each $i=1, \ldots, n-1$, \[S_i \subseteq \{y \in K_2: \abs{Tg(y)-\xi}<\ep\} \cap \{y \in K_2: \abs{Th(y)-\xi}<\ep\},\]
    \item\label{it:almost constant 2}
    \[S_n \subseteq \{y \in K_2: \abs{Tg(y)-\xi}<\ep\}.\]
\end{enumerate}
Then, $\norm{T} \geq n+\sqrt{(n-1)(n+3)}$.
\end{lemma}

\begin{proof}
Put
\[
 \Delta=n+\sqrt{(n-1)(n+3)},\qquad
 D=\norm T,\qquad s=\abs\xi,
\]
and suppose that $D<\Delta$. Choose $\ep_j\downarrow0$ and write the
corresponding objects as $f_i^j,g_j,h_j,U_j$. Set
\[
 b_j=\norm{Tf_n^j},\qquad
 m_j=\max_{i<n}\norm{Tf_i^j}.
\]
Since $\norm{T^{-1}}=1$, all these numbers belong to $[1,D]$.
Passing to a subsequence, we may therefore assume that
$b_j\to b$ and $m_j\to m$.

Fix $A>1$ and put $a=(A-1)/(n-1)$. By \ref{it:supports}, there is a point $x\in K_1$ such that
$f_i^j(x)=1$ for every $i=1,\ldots,n$: indeed, choose $x$ with
$f_1^j(x)=1$ and proceed inductively. Since $0\leq f_i^j\leq1$, it
follows that
\[
\norm{f_n^j+a\sum_{i<n}f_i^j}
=1+a(n-1)=A.
\]
By \ref{it:mapped supports}, for each $y\in K_2$ there are three
possibilities: no index $k$ satisfies
$\abs{Tf_k^j(y)}\geq\ep_j$, the unique such index is $k=n$, or it is
some $k<n$. In the last two cases, we have that $\abs{Tf_k^j(y)}$ is bounded from the above by at most $b_j$ or $m_j$, respectively, while all the remaining
terms have absolute value less than their coefficient times $\ep_j$.
Since the sum of all coefficients is $A$, in all three cases
\[
\abs{Tf_n^j(y)+a\sum_{i<n}Tf_i^j(y)}
\leq\max\{b_j,am_j\}+A\ep_j.
\]
Taking the supremum over $y$ and using $\norm{T^{-1}}=1$, we obtain
\[
A
\leq\norm{T\left(f_n^j+a\sum_{i<n}f_i^j\right)}
\leq\max\{b_j,am_j\}+A\ep_j.
\]
After passing to the limit, we have
\begin{equation}
\label{eq:limit-index-estimate}
 A\leq\max\{b,am\}.
\end{equation}
If $b\leq1$, then in fact $b=1$. Hence, for every $A>1$, we have
$b<A$, and \eqref{eq:limit-index-estimate} forces $am\geq A$.
Consequently,
\[
m\geq\frac{A}{a}=\frac{A(n-1)}{A-1},
\]
which is impossible as $A\downarrow1$, since $m\leq D$. Hence $b>1$. Applying
\eqref{eq:limit-index-estimate} with $A>b$ we have $m\ge b/a$ and then letting $A\downarrow b$
gives
\begin{equation}
\label{eq:short-m-lower}
 m\ge \frac{(n-1)b}{b-1}.
\end{equation}

For each $j$, choose $i_j<n$ such that
$\norm{Tf_{i_j}^j}=m_j$, and choose points $y_j,z_j\in K_2$ satisfying $\abs{Tf_{i_j}^j(y_j)}=m_j$ and $\abs{Tf_n^j(z_j)}=b_j$. For all sufficiently large $j$, we have $y_j\in S_{i_j}$ and
$z_j\in S_n$. Since the sets $S_i$ are pairwise disjoint,
$\abs{Tf_n^j(y_j)}<\ep_j$. Put $q=T^{-1}(\mathbf1_{K_2})$. Since $i_j<n$, conditions
\ref{it:almost 0} and \ref{it:mapped supports}, together with
$Tq=\mathbf{1}_{K_2}$, give
\[
\abs{T(q\chi_{K_1\setminus U_j})(y_j)}
=\abs{1-T(q\chi_{U_j})(y_j)}\stackrel{\ref{it:almost 0}}{>}1-\ep_j.
\]
Choose $\lambda_j\in\{-1,1\}$ so that
$\lambda_jT(q\chi_{K_1\setminus U_j})(y_j)$ and
$-2Tf_{i_j}^j(y_j)$ have the same sign. Then
\[
\begin{aligned}
\norm{T\left(\lambda_jq\chi_{K_1\setminus U_j}
              +f_n^j-2f_{i_j}^j\right)} & \geq
\abs{T(q\chi_{K_1\setminus U_j})(y_j)}
+2\abs{Tf_{i_j}^j(y_j)}
-\abs{Tf_n^j(y_j)}\\
 & \geq1+2m_j-2\ep_j.
\end{aligned}
\]
The function inside $T$ has norm one by
\ref{it:supports} and \ref{it:supports2}; hence, letting $j\to\infty$,
\begin{equation}
\label{eq:first-limit-estimate}
 D\geq1+2m.
\end{equation}
Similarly, condition \ref{it:almost constant}, evaluated at $y_j$ and
with the sign chosen appropriately, gives
\[
 \norm{T\left(\lambda_jh_j+f_n^j-2f_{i_j}^j\right)}
 \geq |Th_j(y_j)| + 2\abs{Tf_{i_j}^j(y_j)}
-\abs{Tf_n^j(y_j)} \geq s+2m_j-2\ep_j.
\]
Again the function inside $T$ has norm one, and therefore
\begin{equation}
\label{eq:second-limit-estimate}
 D\geq s+2m.
\end{equation}
Finally, evaluating $T(g_j-2f_n^j)$ at $z_j$ and using
\ref{it:almost constant 2}, we obtain
\[
 D\geq\norm{T(g_j-2f_n^j)}\geq2b_j-s-\ep_j,
\]
and hence
\begin{equation}
\label{eq:third-limit-estimate}
 D\geq2b-s.
\end{equation}

It remains only to prove the following numerical claim.

\begin{claim}
Let $n\geq2$, $b>1$, $s\geq 0$ and $m,D\in(0,\infty)$. If \eqref{eq:short-m-lower},
\eqref{eq:first-limit-estimate},
\eqref{eq:second-limit-estimate}, and
\eqref{eq:third-limit-estimate} hold, then $D\geq\Delta=n+\sqrt{(n-1)(n+3)}$.
\end{claim}

\begin{proof}[Proof of the claim]
Suppose that $D<\Delta$ and put $c=\max\{1,s\}$. From
\eqref{eq:first-limit-estimate} and
\eqref{eq:second-limit-estimate} we obtain $2m\leq D-c$. Together with \eqref{eq:short-m-lower}, this gives
\begin{equation}\label{eq:short-b-lower}
b\geq\frac{D-c}{D-c-2n+2},
\end{equation}
where the denominator is positive since
\eqref{eq:short-m-lower} implies $m>n-1$, so $D-c\ge 2m > 2(n-1)$.

The roots of
\[
P(t)=-t^2+2nt+2n-3
\]
are $n\pm\sqrt{(n-1)(n+3)}$. The smaller root is negative, while the larger one is $\Delta$. Hence $P(D)>0$.

If $s\leq1$, then $c=1$, and
\eqref{eq:short-b-lower} and \eqref{eq:third-limit-estimate} give
\[
D\geq2b-1
\geq2\frac{D-1}{D-2n+1}-1>D,
\]
because
\[
2\frac{D-1}{D-2n+1}-1-D
=\frac{P(D)}{D-2n+1}>0.
\]

Thus $s>1$, so $c=s$. Put $r=D-s$. Then $2n-2<r<D-1$, and \eqref{eq:short-b-lower} and
\eqref{eq:third-limit-estimate} yield
\[
D\geq \frac{2r}{r-2n+2} - D + r,
\]
which gives
\[
D\geq F(r):=\frac{r}{2}+\frac r{r-2n+2}.
\]
Moreover, $F$ is strictly decreasing on $[r,D-1]$, since, for
$t\in(r,D-1)$,
\[
0<t-2n+2<\Delta - 2n + 1 \le 2\sqrt{n-1}
\quad\text{and hence}\quad
F'(t)=\frac12-\frac{2n-2}{(t-2n+2)^2}<0.
\]
Consequently,
\[
D\geq F(r)>F(D-1)>D,
\]
where
\[
F(D-1)-D
=\frac{P(D)}{2(D-2n+1)}>0.
\]
This contradiction proves the claim.
\end{proof}

The claim contradicts our assumption that $D<\Delta$ and finishes the
proof of the lemma.
\end{proof}

\begin{prop}
\label{prop:infinite splitting old}
Let $K_1, K_2$ be nonempty compact spaces. Suppose that $T:\C(K_1) \rightarrow \C(K_2)$ is an isomorphic embedding, let $L_1 \subseteq K_1$, $U$ be an open set containing $L_1$, and let $L_2 \subseteq K_2$ be a compact set. If there exists an ordinal $\alpha$ such that $\height(L_1)>\omega^{\alpha}$ and $\height(L_2)\leq \omega^{\alpha}$, then for each $\ep>0$ there exist functions $(f_n)_{n \in \en} \subseteq \C(K_1, [0, 1])$, pairwise disjoint open sets $(U_n)_{n \in \en} \subseteq U$, and points $(x_n)_{n \in \en} \subseteq L_1$ such that for each $n \in \en$, $f_n=1$ on an open neighbourhood of $x_n$, $f_n=0$ on $K_1 \setminus U_n$, and $\norm{Tf_n|_{L_2}}<\ep$. 

Moreover, if $K_1$ is zero-dimensional, then the sets $(U_n)_{n \in \en}$ can be taken clopen.
\end{prop}

\begin{proof}

That there exists one such function is the content of \cite[Proposition 2.4(a)]{rondos-somaglia}. To show that there are infinitely many such functions that are disjointly supported, we need to just slightly modify the argument from the proof of \cite[Proposition 2.4]{rondos-somaglia}.

First, we assume that $\alpha=0$. Then, the assumptions of the proposition imply that $L_2$ is finite, while $L_1$ is infinite. Since $L_1$ is infinite, there exist functions $(f_n)_{n \in \en} \subseteq \C(K_1, [0, 1])$, pairwise disjoint open sets $(U_n)_{n \in \en} \subseteq U$, and points $(x_n)_{n \in \en} \subseteq L_1$ such that for each $n \in \en$, $f_n=1$ on an open neighbourhood of $x_n$, $f_n=0$ on $K_1 \setminus U_n$. Let $\ep>0$. It is sufficient to show that any finite subcollection $(f_j)_{j\in M}$ with $|M|>\|T\||L_2|/\ep$  contains $j_0\in M$ such that $\|Tf_{j_0}|_{L_2}\|<\ep$. 

Suppose for a contradiction that $\|Tf_j|_{L_2}\| \geq \epsilon$ for all $j\in M$. Then for each $j$, there exists $y_j \in L_2$ such that $|Tf_j(y_j)| \geq \epsilon$. By the pigeonhole principle, there exists $y_0 \in L_2$ and a subset $M' \subseteq M$ with $|M'| \geq \frac{|M|}{|L_2|}$ such that $y_j = y_0$ for all $j \in M'$. Choosing signs $\sigma_j \in \{\pm 1\}$ such that $\sigma_j Tf_j(y_0) = |Tf_j(y_0)|$, we obtain
\[  |M'|\ep \leq \sum_{j \in M'} |Tf_j(y_0)| = \Big|\sum_{j \in M'} \sigma_j Tf_j(y_0)\Big| \leq \Big\|T\Big(\sum_{j \in M'} \sigma_j f_j\Big)\Big\| \leq \|T\| \Big\|\sum_{j \in M'} \sigma_j f_j\Big\|. \]
Since the functions $f_j$ have disjoint supports, $\left\|\sum_{j \in M'} \sigma_j f_j\right\| = 1$. Thus, $\frac{|M| \ep}{|L_2|}\leq |M'|\ep \leq \|T\|$, which contradicts our choice of $M$.

Next, we assume that $\alpha=\beta+1$ is a successor ordinal. Hence, since $\height(L_2) \leq \omega^{\alpha}$ and $\height(L_2)$ is a successor ordinal, we even have $\height(L_2)<\omega^{\alpha}$. Thus, there exists $k \in \en$ such that $\height(L_2) \leq \omega^{\beta}k$. For a given $\ep>0$ we find a natural number $N \in \en$ such that 
\[\frac{k\norm{T}+(n-k)\frac{\ep}{2}}{n}<\ep\]
for each $n \geq N$.

Next, since $\height(L_1)>\omega^{\alpha}$,  we can find  pairwise disjoint open sets $(U_n)_{n=N}^{\infty} \subseteq U$ such that for each $n \geq N$, $\height(U_n \cap L_1) > \omega^{\beta}n$. Finally, by \cite[Proposition 2.4(b)]{rondos-somaglia}, for each $n \geq N$, we can find a point $x_n \in U_n \cap L_1$ and a function $f_n \in \C(K_1, [0, 1])$ such that $f_n=1$ on a neighbourhood of $x_n$, $f_n=0$ on $K_1 \setminus U_n$ and 
\[\norm{Tf_n|_{L_2}} \leq \frac{k\norm{T}+(n-k)\frac{\ep}{2}}{n}<\ep.\]
This finishes the proof in the successor case. 

Further, the case when $\alpha$ is a limit ordinal follows trivially from the successor case. Indeed, if $\alpha$ is a limit ordinal and $\height(L_2) \leq \omega^{\alpha}$, then, since $\height(L_2)$ is a successor ordinal, $\height(L_2) <\omega^{\alpha}$, hence, there exists a successor ordinal $\beta<\alpha$ such that $\height(L_2) \leq \omega^{\beta}$. At the same time, $\height(L_1)>\omega^{\beta}$, hence the conclusion of the lemma follows from what was proved above. 

Finally, it is easy to see that in the argument above, the sets $(U_n)_{n \in \en}$ can be taken clopen if $K_1$ is zero-dimensional.

\end{proof}

Before proceeding further, we now recall that for any compact space $K$ and a point $x \in K$, the characteristic function $\chi_x$ can be naturally viewed as an element of $\C(K)^{**}$ via the formula 
$\chi_x(\mu)=\mu(\{x\})$, $\mu \in M(K)$. Furthermore, it is a standard fact that the function $\chi_x$ is a weak$^*$ limit of the net of functions $\{f_U: U \text{ is a neighbourhood of } x\}$, where the function $f_U \in \C(K, [0, 1])$ vanishes outside of the set $U$, and satisfies $f_U=1$ on some open neighbourhood of $x$. Hence, if $T: \C(K_1) \rightarrow \C(K_2)$ is an isomorphism, then $\{Tf_U\}$ converges weak$^*$ to $T^{**}\chi_x$, by the weak$^*$-weak$^*$ continuity of $T^{**}$.

\begin{thm}
\label{thm:lower estimate general isomorphism}
Let $K_1, K_2$ be scattered compact spaces, and assume that for some ordinal $\alpha$ and $n \in \en, n \geq 2$, $\height(K_1)>\omega^{\alpha}\cdot n$, and $K_2^{(\omega^{\alpha})}=\{y\}$ for some $y \in K_2$. Then
\[d_{BM}(\C(K_1), \C(K_2)) \geq n+\sqrt{(n-1)(n+3)}.\]
\end{thm}

\begin{proof}
We fix an isomorphism $T:\C(K_1) \rightarrow \C(K_2)$ with $\norm{T^{-1}}=1$, and let $\ep>0$ be arbitrary. To prove the required lower bound, it is enough to find the objects from the statement of Lemma \ref{lemma:the lower estimate C(n)}.

To this end, we pick a point $x \in K_1^{(\omega^{\alpha}\cdot n)}$, and we denote $\xi=T^{**}(\chi_{\{x\}})(\delta_y)$, where $\delta_y$ stands for the Dirac measure at the point $y$. Further, we find an open neighbourhood $V$ of $x$ and a function $g \in \C(K_1, [0, 1])$ such that $g=1$ on $V$, and such that $\abs{Tg(y)-\xi}<\ep$. 

Further, we denote $L_1=V^{(\omega^{\alpha}(n-1))}$ and $L_2=\{z \in K_2: \abs{Tg(z)-\xi} \geq \ep\}$. Then, since $x \in L_1$, $\height(L_1)>\omega^{\alpha}$, and, since $y \notin L_2$, $\height(L_2) \leq \omega^{\alpha}$. Notice also that the set $L_2$ is compact. Consequently, by Proposition \ref{prop:infinite splitting old} we can find pairwise disjoint clopen sets $(U_k)_{k \in \en} \subseteq V$ (notice that $K_1$ is zero-dimensional, since it is scattered), functions $g_k \in \C(K_1, [0, 1])$, and points $(x_k)_{k \in \en} \subseteq L_1$ such that for each $k \in \en$, $g_k=1$ on an open neighbourhood $W_k$ of $x_k$, $g_k=0$ on $K_1 \setminus U_k$, and $\norm{Tg_k|_{L_2}}<\ep$. Since the sets $U_k$ are pairwise disjoint, we may assume that none of them contains the point $x$. Next, by a standard pigeonhole-like argument,   there exists $k \in \en$ such that
\[\abs{T(T^{-1}({\mathbf 1}_{K_2})\cdot \chi_{U_k})(y)}< \ep \text{ and } \abs{Tg_k(y)}<\ep.\] 
We denote $U=U_k$, $V_n=W_k$ and $f_n=g_k$. Thus, so far, we have found the functions $f_n$ and $g$, and the clopen set $U$ from the statement of Lemma \ref{lemma:the lower estimate C(n)}, and, we
have the conditions \ref{it:supports2} and \ref{it:almost constant 2} satisfied.

In the next step, we find a function $h \in \C(K_1, [0, 1])$ such that $\abs{Th(y)-\xi}<\ep$ and $\spt h \cap U=\emptyset$, so, also the condition \ref{it:supports3} is satisfied.

Now, we aim to find the required functions $f_1, \ldots, f_{n-1} \in \C(K_1, [0, 1])$ by reverse induction. Thus, we assume that for some $0 \leq k<n-1$, we have found the functions $f_{n-k}, \ldots, f_n$, and open sets $V_{n-k}, \ldots, V_n$ such that for each $i=n-k, \ldots, n$, we have $\height(V_i)>\omega^{\alpha}(i-1)$, $f_i=1$ on $V_i$, the sets \[S_i=\{z \in K_2: \abs{Tf_i(z)} \geq \ep\}_{i=n-k}^n\] are pairwise disjoint and do not contain the point $y$, and such that   
\[S_i \subseteq \{z \in K_2: \abs{Tg(z)-\xi} < \ep\}, \quad i=n-k, \ldots, n,\] 
\[S_i \subseteq \{z \in K_2: \abs{T(T^{-1}({\mathbf 1}_{K_2})\cdot \chi_{U})(z)} < \ep\}, \quad i=n-k, \ldots, n-1,\]and 
\[S_i \subseteq \{z \in K_2: \abs{Th(z)-\xi} < \ep\}, \quad i=n-k, \ldots, n-1.\]

Now, to find the function $f_{n-k-1}$, we denote $L_1=V_{n-k}^{(\omega^{\alpha}(n-k-2))}$, and
\[
\begin{aligned}
L_2=&\{z \in K_2: \abs{Tg(z)-\xi} \geq \ep \text{ or } \abs{Th(z)-\xi} \geq \ep \text{ or } \abs{T(T^{-1}({\mathbf 1}_{K_2})\cdot \chi_U)(z)} \geq \ep \text{ or } &\\& \abs{Tf_i(z)} \geq \ep \text{ for some } i=n-k, \ldots, n\}.
\end{aligned}
\]

Then, since $\height(V_{n-k})>\omega^{\alpha}(n-k-1)$, $\height(L_1) > \omega^{\alpha}$, while $\height(L_2)\leq\omega^{\alpha}$, since $y \notin L_2$. And, the set $L_2$ is again compact. Consequently, using Proposition \ref{prop:infinite splitting old} and using pigeonhole-like argument similarly as above, we
can find a function $f_{n-k-1} \in \C(K_1, [0, 1])$ such that $f_{n-k-1}$ vanishes outside of the set $V_{n-k}$, $f_{n-k-1}=1$ on an open neighbourhood $V_{n-k-1}$ of some point from the set $V_{n-k}^{(\omega^{\alpha}(n-k-2))}$ (which, in particular, implies that $\height(V_{n-k-1})>\omega^{\alpha}(n-k-2)$), and such that $\norm{Tf_{n-k-1}|_{L_2}}< \ep$   and $|Tf_{n-k-1}(y)|<\varepsilon$. This finishes the induction step.

Now, from the construction it follows that also all the remaining required conditions \ref{it:supports}, \ref{it:mapped supports}, \ref{it:almost 0}, and \ref{it:almost constant} from Lemma \ref{lemma:the lower estimate C(n)} are satisfied, and hence the proof is finished.
\end{proof}

\begin{cor}
\label{cor:lower-estimate-special-ordinals}
For each ordinal $\alpha$ and $n \in \en$, 
\[
d_{BM}\bigl(\C(\omega^{\omega^{\alpha}}),
\C(\omega^{\omega^{\alpha}\cdot n})\bigr)
\geq n+\sqrt{(n-1)(n+3)}.
\]
\end{cor}
\begin{proof}
Apply Theorem \ref{thm:lower estimate general isomorphism} on $K_1=[1, \omega^{\omega^{\alpha}\cdot n}]$ and $K_2=[1, \omega^{\omega^{\alpha}}]$ (we have $\height(K_1)=\omega^{\alpha}\cdot n+1$ and $K_2^{(\omega^{\alpha})}=\{\omega^{\omega^{\alpha}}\}$).    
\end{proof}

\begin{proof}[Proof of Theorem \ref{thm:the exact values of the distances}.]
The result follows immediately from Corollary~\ref{cor:lower-estimate-special-ordinals} combined with Proposition~\ref{prop:constructing-isomorphism-new-attempt}(ii).
\end{proof}

\subsection{A general lower bound for the positive distance}
In this subsection we prove Proposition~\ref{prop:lower estimate same heights}, a lower estimate for compact spaces whose last nonempty Cantor--Bendixson derivatives have different finite cardinalities. We start with Lemma~\ref{lem:wlogOneToOne}, which implies that when dealing with the one-sided positive Banach--Mazur distance, we may without loss of generality assume that the constant-one function is fixed by our isomorphism.

\begin{lemma}\label{lem:wlogOneToOne}Let $K_1,K_2$ be compact spaces and let $T:C(K_1)\to C(K_2)$ be a positive surjective isomorphism. Then $\norm{T} = \norm{T(\boldsymbol{1})}$ and $T(\boldsymbol{1})\geq \frac{1}{\norm{T^{-1}}}\cdot \boldsymbol{1}$.
Moreover, there exists a positive surjective isomorphism $S:C(K_1)\to C(K_2)$ such that $\norm{S}=1$, $\norm{S^{-1}}\leq \norm{T}\cdot \norm{T^{-1}}$ and $S(\boldsymbol{1})=\boldsymbol{1}$.
\end{lemma}
\begin{proof}
We may without loss of generality assume that $\norm{T}=1$. Now, since $T(\boldsymbol{1})\geq 0$ and $\norm{T(\boldsymbol{1})}\leq 1$, we have $0\leq T(\boldsymbol{1})\leq \boldsymbol{1}$. Furthermore, for any $f\in B_{C(K_1)}$ we have $0\leq T(f^+)\leq T(\boldsymbol{1})$ and $0\leq T(f^-)\leq T(\boldsymbol{1})$, which implies
\[
-T(\boldsymbol{1})(y)\leq T(f^+)(y)-T(f^-)(y) = Tf(y)\leq T(\boldsymbol{1})(y),\quad y\in K_2,
\]
so we obtain $|Tf|\leq T(\boldsymbol{1})$ (in particular, we have $\norm{Tf} \leq \norm{T(\boldsymbol{1})}$ and since $f\in B_{C(K_1)}$ was arbitrary, we obtain $\norm{T(\boldsymbol{1})} = \norm{T}=1$). Now, given $y\in K_2$ we pick $g\in S_{C(K_2)}$ with $g(y)=1$ and then for $f:=\frac{T^{-1}g}{\norm{T^{-1}}}\in B_{C(K_1)}$ we obtain $Tf(y) = \frac{1}{\norm{T^{-1}}}$ and therefore $T(\boldsymbol{1})(y) \geq Tf(y) = \frac{1}{\norm{T^{-1}}}$.

For the ``Moreover'' part it suffices to put $Sf:=\frac{Tf}{T(\boldsymbol{1})}$. Then it is easy to check that $S$ is positive surjective isomorphism with $S(\boldsymbol{1})=\boldsymbol{1}$ and by the already proven part we have $\norm{S}=\norm{S(\boldsymbol{1})} = 1$. For $g\in C(K_2)$ we have
\[
\|S^{-1}g\|= \|T^{-1}(T(\boldsymbol{1})\cdot g)\|\le \|T^{-1}\|\cdot \|T(\boldsymbol{1})\cdot g\|
\le \|T^{-1}\|\cdot \|T(\boldsymbol{1})\|\cdot \|g\|=\|T^{-1}\|\cdot \|T\|\cdot \|g\|.
\]
This yields $\|S^{-1}\|\le \|T\|\cdot \|T^{-1}\|$.
\end{proof}

Our final argument is then based on the following easy observation.

\begin{lemma}\label{lem:triv1}Let $T:\C(K)\to \C(L)$ be a surjective positive isomorphism with $T\boldsymbol{1} = \norm{T}\cdot \boldsymbol{1}$ and $\alpha > 0$. Assume there exists a nonnegative function $f\in S_{\C(K)}$ with $\norm{Tf}\leq \alpha\norm{T}$. Then $\norm{T}\cdot \norm{T^{-1}}\geq \tfrac{2}{\alpha}-1$.
\end{lemma}
\begin{proof}Assume $0\leq f\leq 1$ with $\norm{f}=1$ is such that $\norm{Tf}\leq \alpha\norm{T}$. Since $0\leq \tfrac{Tf}{\norm{Tf}}\leq 1$, we have
\[
\norm{T^{-1}}\geq \norm{T^{-1}(1-2\tfrac{Tf}{\norm{Tf}})} = \norm{\tfrac{1}{\norm{T}}-\tfrac{2}{\norm{Tf}}f} \geq \tfrac{1}{\norm{T}}(\tfrac{2\norm{T}}{\norm{Tf}}-1)\geq \tfrac{1}{\norm{T}}(\tfrac{2}{\alpha}-1).
\]
\end{proof}

\begin{prop}
\label{prop:lower estimate same heights}    
Assume that $K_1, K_2$ are compact spaces such that for some ordinal $\alpha$ and $n \in \en$, $\abs{K_1^{(\alpha)}}=n$ and $\abs{K_2^{(\alpha)}}=1$. Then, \[d_{BM}^+(C(K_1),C(K_2))\geq 2n-1.\]
\end{prop}

\begin{proof}
For $n=1$ this is trivial, so we may assume $n\geq 2$. Pick a positive isomorphism $T:C(K_1)\to C(K_2)$. By Lemma~\ref{lem:wlogOneToOne}, without loss of generality, we may assume $T\boldsymbol{1}= \boldsymbol{1}$ and $\norm{T}=1$. Next, since $K_1$ is scattered, it is zero-dimensional, hence we find pairwise disjoint clopen sets $\{U_i\}_{i=1}^n$ of height $\alpha+1$, covering $K_1$. We denote $f_i=\chi_{U_i}$ for $i\le n$, and let $y_0 \in K_2$ denote the point satisfying $K_2^{(\alpha)}=\{y_0\}$. Since \[0\leq T(\sum_{i=1}^n f_i)(y_0)=T(1)(y_0) = 1,\] 
there is $i\in \{1,\ldots,n\}$ such that $T(f_i)(y_0)\leq \tfrac{1}{n}$. Further, we fix
$0<\varepsilon<\min\{\tfrac{1}{n},\tfrac{1}{\norm{T^{-1}}}\}$, and we consider the compact set 
\[L=\{y \in K_2\setsep |Tf_i(y)-Tf_i(y_0)|\geq \varepsilon\}.\] Then, since $y_0 \notin L$, we have $\height(L)<\alpha+1=\height(U_i)$. Consequently, by Proposition \ref{prop:basic properties of positive embeddings}(a), we can find a norm-$1$ function $g \in \C(K_1, [0, 1])$ supported by $U_i$, such that $\norm{Tg|_L}<\varepsilon$. We also have
\[
\norm{Tg|_{L^c}}\leq \norm{Tf_i|_{L^c}}\leq T(f_i)(y_0) + \varepsilon\leq \tfrac{1}{n} + \varepsilon,
\]
and so $\norm{Tg}\leq \alpha(\varepsilon)\norm{T}$ for $\alpha(\varepsilon) = \tfrac{1}{n} + \varepsilon$. By Lemma~\ref{lem:triv1}, this implies that $\norm{T}\cdot\norm{T^{-1}}\geq \tfrac{2}{\alpha(\varepsilon)} - 1$ and since $\varepsilon>0$ could be chosen arbitrarily small, we obtain $\norm{T}\cdot\norm{T^{-1}}\geq \lim_{\varepsilon\to 0^+}\tfrac{2}{\alpha(\varepsilon)} - 1 = 2n-1$.
\end{proof}

\subsection{A two-branch upper construction}

\begin{proof}[Proof of Proposition~\ref{prop:pd-higher-two-upper}]
Apply Lemma~\ref{lemma:split into subintervals} with
$\Gamma=\{1,2\}$, and let $I_1,I_2$ be the resulting partition of
$[0,\omega^\alpha)$. Put
\[
x_0=\omega^\alpha,\qquad
y_1=\omega^\alpha,\qquad
y_2=\omega^\alpha\cdot2,
\]
and let
\[
J_1=[0,\omega^\alpha),\qquad
J_2=[\omega^\alpha+1,\omega^\alpha\cdot2).
\]
For $i=1,2$, choose homeomorphisms
\[
r_i:[0,\omega^\alpha]\longrightarrow\overline{I_i},
\qquad
s_i:[0,\omega^\alpha]\longrightarrow\overline{J_i}
\]
such that $r_i(\omega^\alpha)=x_0$ and
$s_i(\omega^\alpha)=y_i$. Such homeomorphisms exist because
$\overline{I_i}$ and $\overline{J_i}$ are the one-point compactifications of
spaces homeomorphic to $[0,\omega^\alpha)$. Put
\[
\bar p_i:=s_i\circ r_i^{-1}:\overline{I_i}\longrightarrow
\overline{J_i},
\qquad
p_i:=\bar p_i|_{I_i}.
\]
Then $p_i:I_i\to J_i$ is a homeomorphism and
$I_i\ni\xi\to x_0$ if and only if $p_i(\xi)\to y_i$.

Let
\[
K_0=[0,\omega^\alpha],
\qquad L_0=[0,\omega^\alpha\cdot2],
\]
and adjoin two isolated points $a_1,a_2$ to $K_0$ and one isolated point $z$
to $L_0$. By the Mazurkiewicz--Sierpi\'nski classification,
\[
\widehat K=K_0\sqcup\{a_1,a_2\}\simeq[0,\omega^\alpha],
\qquad
\widehat L=L_0\sqcup\{z\}\simeq[0,\omega^\alpha\cdot2].
\]

Fix $0<\lambda<1$. For $f\in\C(\widehat K)$ define
$Tf\in\C(\widehat L)$ by
\[
\begin{split}
Tf(y_i)&=\lambda f(x_0)+(1-\lambda)f(a_i)
        \qquad(i=1,2),\\
Tf(z)&=\tfrac12\bigl(f(a_1)+f(a_2)\bigr),\\
Tf(p_i(\xi))&=(1-\lambda)f(a_i)+\lambda f(\xi)
        \qquad(\xi\in I_i,\ i=1,2).
\end{split}
\]
Since $\bar p_i(x_0)=y_i$, the formula on $J_i$ extends continuously to
$y_i$ with the prescribed value. Hence $Tf$ is continuous. Moreover, $T$ is
positive, $T\boldsymbol{1}=\boldsymbol{1}$, and $\norm{T}=1$.

For $g\in\C(\widehat L)$ define $Sg\in\C(\widehat K)$ by
\[
\begin{split}
Sg(x_0)&=\frac{g(y_1)+g(y_2)}{2\lambda}
         -\frac{1-\lambda}{\lambda}g(z),\\
Sg(a_1)&=\frac{g(y_1)-g(y_2)}{2(1-\lambda)}+g(z),\\
Sg(a_2)&=\frac{-g(y_1)+g(y_2)}{2(1-\lambda)}+g(z),\\
Sg(\xi)&=\frac{-g(y_1)+g(y_2)-2(1-\lambda)g(z)
                   +2g(p_1(\xi))}{2\lambda}
                   \qquad(\xi\in I_1),\\
Sg(\xi)&=\frac{g(y_1)-g(y_2)-2(1-\lambda)g(z)
                   +2g(p_2(\xi))}{2\lambda}
                   \qquad(\xi\in I_2).
\end{split}
\]
Indeed, as $\xi\to x_0$ in either $I_1$ or $I_2$, the corresponding
formula tends to
\[
\frac{g(y_1)+g(y_2)-2(1-\lambda)g(z)}{2\lambda}=Sg(x_0).
\]
Thus $Sg$ is continuous. Direct substitution gives $ST=\Id$ and $TS=\Id$,
and
\[
\norm{S}
=\max\left\{1+\frac1{1-\lambda},\frac3\lambda-1\right\}.
\]
For $\lambda=(3-\sqrt3)/2$, both terms in the maximum are equal to
$2+\sqrt3$. Consequently,
\[
\norm{T}\norm{T^{-1}}\leq2+\sqrt3,
\]
which proves the assertion.
\end{proof}

The next section proves that the upper bound in
Proposition~\ref{prop:Construction with better bound} is optimal for every
nonzero countable exponent. It also determines the opposite-orientation
distances at rank one and at principal higher ranks.

\section{Exact positive distances for finite ordinal multiples}\label{sec:positive-distances}

Let us recall the upper bounds we proved so far. For every nonzero countable ordinal $\alpha$ and every integer $k\geq 2$, we have by Proposition~\ref{prop:Construction with better bound} and Proposition~\ref{prop:constructing-isomorphism-new-attempt} applied to $\beta = k$ and $n=1$
\[
d_{BM}^+(\C(\omega^\alpha\cdot k),\C(\omega^\alpha))
\leq k+\sqrt{(k-1)(k+3)}\quad \text{and}\quad d_{BM}^+(\C(\omega^\alpha),\C(\omega^\alpha\cdot k))\leq 2+\sqrt5,
\]
and by Proposition~\ref{prop:pd-higher-two-upper} we have  
\[
d_{BM}^+(\C(\omega^\alpha),\C(\omega^\alpha\cdot2))\leq 2+\sqrt3.
\]
Hence, in order to prove Theorem~\ref{thm:Intro3} it suffices to prove the corresponding lower bounds, which is the aim of this entire section.

\subsection{Normalization and localization}
\label{subsec:pd-localization-principles}

Recall that an operator $T:\C(K)\to\C(L)$ is called \emph{unital} if
$T\mathbf 1_K=\mathbf 1_L$. By Lemma~\ref{lem:wlogOneToOne}, when
computing the positive Banach--Mazur distance we may restrict ourselves
to positive unital isomorphisms. For such an isomorphism $T$ we have,
again by Lemma~\ref{lem:wlogOneToOne}, $\norm T=1$.

We shall repeatedly use the following elementary consequences of
unitality. For every $y\in L$, the measure
$\mu_y:=T^*\delta_y$ is positive and $\mu_y(K)=T\mathbf 1_K(y)=1$, so $\mu_y$ is a probability measure. Moreover, if $S=T^{-1}$, then $S\mathbf 1_L=\mathbf 1_K$, and hence $(S^*\delta_x)(L)=1$ for every $x\in K$.

We shall also use without further comment that, for a countable compact
space $X$, the identification
\[
M(X)=\ell_1(X),\qquad
\lambda\longleftrightarrow (\lambda(x))_{x\in X},\qquad
\lambda(x)=\lambda(\{x\}),
\]
is isometric. Thus $\lambda=\sum_{x\in X}\lambda(x)\delta_x$ and
$\norm\lambda=\sum_{x\in X}\abs{\lambda(x)}$.

\begin{lemma}[Shrinking-neighbourhood escape]
\label{lem:pd-shrinking-escape}
Let $X$ be compact, let $y\in X$, and let $(U_n)$ be a decreasing compact
neighbourhood basis at $y$. If $\Gamma_n\xrightarrow{w^*}\Gamma$ in $M(X)$
and $\liminf_n\Gamma_n^+(U_n)\ge a\ge0$, then
$\liminf_n\norm{\Gamma_n}\ge a+\norm{\Gamma-a\delta_y}$.
\end{lemma}

\begin{proof}
Put $L=\liminf_n\norm{\Gamma_n}$. After passing to a subsequence, we may
assume that $\norm{\Gamma_n}\to L$ and
$b_n:=\Gamma_n^+(U_n)\to b\ge a$. Put
$\theta_n:=\Gamma_n^+\mathbin{\upharpoonright}U_n$. For every $f\in\C(X)$,
\[
\abs{\theta_n(f)-bf(y)}
\le b_n\sup_{x\in U_n}\abs{f(x)-f(y)}
   +\abs{b_n-b}\abs{f(y)}\longrightarrow0.
\]
Thus $\theta_n\xrightarrow{w^*}b\delta_y$. Since
$0\le\theta_n\le\Gamma_n^+$,
$\Gamma_n-\theta_n=(\Gamma_n^+-\theta_n)-\Gamma_n^-$ is the Jordan
decomposition of $\Gamma_n-\theta_n$. Hence
$\norm{\Gamma_n-\theta_n}=\norm{\Gamma_n}-b_n$. By weak-star lower
semicontinuity,
\[
\norm{\Gamma-b\delta_y}
\le\liminf_n\norm{\Gamma_n-\theta_n}
=\lim_n\bigl(\norm{\Gamma_n}-b_n\bigr)=L-b.
\]
Consequently, using
$\norm{\Gamma-a\delta_y}\le\norm{\Gamma-b\delta_y}+b-a$, we obtain
\[
L\ge b+\norm{\Gamma-b\delta_y}
\ge a+\norm{\Gamma-a\delta_y}.\qedhere
\]
\end{proof}

\subsection{Positive maps from $\C(\omega^\alpha\cdot k)$ to
$\C(\omega^\alpha)$}
\label{sec:pd-forward-all-ranks}
The aim of this subsection is to prove the following result.

\begin{thm}\label{thm:pd-forward-all-ranks}
Let $0<\alpha<\omega_1$ and let $k\ge2$ be an integer. Then
\[
d_{BM}^+\bigl(\C(\omega^\alpha\cdot k),\C(\omega^\alpha)\bigr)
=k+\sqrt{(k-1)(k+3)}.
\]
\end{thm}

By Proposition~\ref{prop:Construction with better bound}, it suffices to prove the lower bound, that is, the inequality ``$\geq$''. Before starting with our argument we establish some notation. Fix $\alpha$ and $k$ as above, and put $K=[1,\omega^\alpha\cdot k]$, $L=[1,\omega^\alpha]$,
$d_j=\omega^\alpha\cdot j$ for $1\le j\le k$, and $e=\omega^\alpha\in L$.
Let
\[
h=(k-1)+\sqrt{(k-1)(k+3)},
\qquad
\Delta=h+1=k+\sqrt{(k-1)(k+3)}.
\]
We shall repeatedly use
\begin{equation}\label{eq:pd-h-relation}
h^2=2(k-1)h+4(k-1).
\end{equation}

By Lemma~\ref{lem:wlogOneToOne}, it suffices to consider a positive unital
isomorphism $T:\C(K)\to\C(L)$. Set $S=T^{-1}$ and $t=\norm S$, and write
$\mu_x=T^*\delta_x\in M(K)^+$ and $\Lambda_z=S^*\delta_z\in M(L)$.
The measures $\mu_x$ are probabilities. Decompose
\[
\mu:=\mu_e=\sum_{j=1}^k A_j\delta_{d_j}+\nu,
\qquad
\nu(d_j)=0\quad(1\le j\le k),
\]
and put $R=\norm\nu$. Then
\begin{equation}\label{eq:pd-forward-mass}
\sum_{j=1}^kA_j+R=1.
\end{equation}

\begin{lemma}[Endpoint localization]\label{lem:pd-forward-localization}
For every $1\le j\le k$,
\[
A_j\ge\frac1t,
\qquad
t\ge\frac1{A_j}
 +\norm{\Lambda_{d_j}-\frac1{A_j}\delta_e}.
\]
\end{lemma}

\begin{proof}
Fix $j$ and abbreviate $A=A_j$ and $d=d_j$. Choose decreasing clopen
neighbourhood bases $(V_n)$ at $d$ and $(U_n)$ at $e$, and numbers
$\ep_n\downarrow0$, with $\ep_n<1/t$, so that
\[
\mu(V_n)\le A+\ep_n,
\qquad
\mu_x(V_n)\le A+2\ep_n\quad(x\in U_n).
\]
This is possible because $x\mapsto\mu_x(V_n)=T\chi_{V_n}(x)$ is continuous.
For $F_n=L\setminus U_n$ we have
$\height(F_n)<\alpha+1=\height(V_n)$. Hence
Proposition~\ref{prop:basic properties of positive embeddings}(a) gives
$f_n\in\C(K,[0,1])$ such that
\[
\norm{f_n}=1,
\qquad
\spt f_n\subset V_n,
\qquad
\norm{Tf_n|_{F_n}}<\ep_n.
\]
Choose $z_n\in V_n$ with $f_n(z_n)=1$ and set
$\Gamma_n=\Lambda_{z_n}$. Since $\norm{\Gamma_n}\le t$,
\[\begin{split}
1 & =\Gamma_n(Tf_n)\le \Gamma_n^+(Tf_n)
 \le \|Tf_n|_{U_n}\|\cdot\Gamma_n^+(U_n) + \|Tf_n|_{F_n}\|\cdot\Gamma_n^+(F_n)\\
 & \le \|T\chi_{V_n}|_{U_n}\|\cdot\Gamma_n^+(U_n) + \ep_nt\le(A+2\ep_n)\Gamma_n^+(U_n)+\ep_nt.
\end{split}\]
Thus $A>0$ and $\liminf_n\Gamma_n^+(U_n)\ge1/A$, which already yields
$A\ge1/t$. Since $z_n\to d$, we have $\delta_{z_n}\xrightarrow{w^*}\delta_d$ and hence
$\Gamma_n=S^*\delta_{z_n}\xrightarrow{w^*}\Lambda_d$. Therefore
Lemma~\ref{lem:pd-shrinking-escape} gives the second assertion.
\end{proof}

Set $\Lambda_j=\Lambda_{d_j}$ and $u_j=A_j\Lambda_j(e)$, and define
\[
G(s)=
\begin{cases}
\dfrac{2(1-s)}h,&s\le0,\\[1em]
\dfrac{2+\dfrac4h s(1-s)}{h+2s},&0\le s\le1,\\[1em]
\dfrac{2s}{h+2},&s\ge1.
\end{cases}
\]

\begin{lemma}[Endpoint estimate]\label{lem:pd-forward-endpoint}
If $t<\Delta$, then $A_j>G(u_j)$ for $1\le j\le k$.
\end{lemma}

\begin{proof}
Fix $j$ and write $A=A_j$, $u=u_j$, and $\Lambda=\Lambda_j$. By
Lemma~\ref{lem:pd-forward-localization},
\[
t\ge\frac1A+\norm{\Lambda-\frac1A\delta_e}.
\]
Decompose $\Lambda=(u/A)\delta_e+\eta$, where $\eta(e)=0$. Since
$M(L)=\ell_1(L)$,
$\norm{\Lambda-A^{-1}\delta_e}=|1-u|/A+\norm\eta$. Moreover,
$T^*\eta=\delta_{d_j}-(u/A)\mu$ and $\norm T=1$, so
\begin{equation}\label{eq:pd-forward-basic-u}
t\ge\frac1A+\frac{|1-u|}{A}+\norm\eta
 \ge\frac1A+\frac{|1-u|}{A}
 +\norm{\delta_{d_j}-\frac uA\mu}.
\end{equation}

If $u\le0$, the last measure is positive and has norm $1-u/A$; hence
$t\ge1+2(1-u)/A$. Since $t<h+1$, this gives $A>2(1-u)/h=G(u)$.

If $u\ge1$, then
\[
\delta_{d_j}-\frac uA\mu
=-(u-1)\delta_{d_j}-\frac uA(\mu-A\delta_{d_j})\le0,
\]
because $\mu-A\delta_{d_j}$ is a positive measure. Hence its norm is the
negative of its total mass, and therefore $\norm{\delta_{d_j}-\frac uA\mu}
=-1+\frac{u}{A}$. Using $|1-u|=u-1$ in \eqref{eq:pd-forward-basic-u}, we obtain
$t\ge-1+2u/A$. Since $t<h+1$, it follows that
    $A>2u/(h+2)=G(u)$.

Assume $0<u<1$ and write $A\eta=\sigma^+-\sigma^-$ for its coordinatewise
positive and negative parts. Since $\eta(e)=0$, both $\sigma^+$ and $\sigma^-$ vanish at $e$. Put $P=T^*\sigma^+$ and $N=T^*\sigma^-$. Then $P,N\ge0$ and
\[
P-N=AT^*\eta=A\delta_{d_j}-u\mu
=A(1-u)\delta_{d_j}-u(\mu-A\delta_{d_j}).
\]
The two positive measures in the last difference have disjoint supports.
Hence, coordinatewise in $\ell_1(K)$, there is a positive measure $\rho$ such that
\[
P=A(1-u)\delta_{d_j}+\rho,
\qquad
N=u(\mu-A\delta_{d_j})+\rho.
\]
Writing $\Theta=2\norm\rho$ and using $\norm{T^*}=1$, we obtain
\[
A\norm\eta\ge\norm P+\norm N
=A(1-u)+u(1-A)+\Theta.
\]
Together with \eqref{eq:pd-forward-basic-u}, this gives
\begin{equation}\label{eq:pd-forward-moving-estimate}
t\ge\frac1A+\frac{1-u}{A}
 +\frac{A(1-u)+u(1-A)+\Theta}{A} = \frac{2+A}{A}-2u+\frac{\Theta}{A}.
\end{equation}

Now put
$\theta=\delta_e-\sigma^+/(1-u)-\sigma^-/u$. Since
$\delta_e,\sigma^+,\sigma^-$ have pairwise disjoint supports,
\[
\norm\theta
=1+\frac{\norm{\sigma^+}}{1-u}+\frac{\norm{\sigma^-}}u
\ge1+\frac{\norm P}{1-u}+\frac{\norm N}{u}
=2+\frac{\Theta}{2u(1-u)}.
\]
On the other hand,
\[
T^*\theta
=\mu-\frac{P}{1-u}-\frac{N}{u}=-\left(\frac1{1-u}+\frac1u\right)\rho,
\]
so $\norm{T^*\theta}=\Theta/(2u(1-u))$. Thus,
\[
2+\frac{\Theta}{2u(1-u)}\le\norm{\theta}
\le t\norm{T^*\theta}
=\frac{t\Theta}{2u(1-u)}.
\]
This first implies $\Theta>0$. Moreover,
\[
2\le(t-1)\frac{\Theta}{2u(1-u)},
\]
so $t>1$ and
\[
\Theta\ge\frac{4u(1-u)}{t-1}>\frac{4u(1-u)}h.
\]
Rearranging \eqref{eq:pd-forward-moving-estimate} and, using $A>0$, $t<h+1$ and the lower bound
for $\Theta$, we obtain
\[
A(h+2u)>A(t-1+2u)\ge2+\Theta
>2+\frac4h u(1-u).
\]
Therefore, dividing by $(h+2u)$ we obtain $
A>G(u)$.
\end{proof}

\begin{lemma}[Residual estimate]\label{lem:pd-forward-residual}
If $t<\Delta$, then
\[
R\ge\frac{2(\sum_{j=1}^ku_j-1)_+}{h}.
\]
\end{lemma}

\begin{proof}
The assertion is immediate if $\sum_j u_j\le1$. Assume that
$b:=\sum_j u_j-1>0$. By the identification $M(K)=\ell_1(K)$,
we have $\nu=\sum_{x\in K}\nu(x)\delta_x$ and
$R=\sum_{x\in K}\nu(x)$. Applying $S^*$ to the decomposition of $\mu$,
using $S^*T^*=\Id_{M(L)}$, and evaluating at $e$, we obtain
\[
1=(S^*\mu)(e)
=\sum_{j=1}^ku_j+\sum_{x\in K}\nu(x)\Lambda_x(e).
\]
Consequently, $
\sum_{x\in K}\nu(x)\bigl(-\Lambda_x(e)\bigr)=b > 0$. Thus, $R\neq 0$ and so $R>0$. Dividing the preceding identity by $R$
gives
\[
\frac bR
=\sum_{x\in K}\frac{\nu(x)}R\bigl(-\Lambda_x(e)\bigr),
\qquad
\sum_{x\in K}\frac{\nu(x)}R=1.
\]
Thus $b/R$ is a weighted average of the numbers $-\Lambda_x(e)$, so for
some $x\in K$ with $\nu(x)>0$ we have
$-\Lambda_x(e)\ge b/R$. Put $r=\Lambda_x(e)$. Then $-r\ge\frac bR>0$, so in particular $r<0$. Write
$\Lambda_x=r\delta_e+\eta$, where $\eta(e)=0$. Then
$T^*\eta=\delta_x-r\mu$. Since $r<0$, this is a positive measure of norm
$1-r$, and hence $\norm\eta\ge1-r$. Moreover, the $\ell_1$
decomposition at $e$ gives $\norm{\Lambda_x}=-r+\norm\eta$. Therefore
\[
t\ge\norm{\Lambda_x}
\ge-r+1-r
=1-2r
\ge1+\frac{2b}{R}.
\]
Since $t<h+1$, it follows that $R>2b/h$, as required.
\end{proof}

It remains to combine the endpoint and residual estimates. This is the following purely numerical fact.

\begin{lemma}[Numerical estimate]\label{lem:pd-forward-numerical}
For all $u_1,\ldots,u_k\in\er$,
\[
\sum_{j=1}^kG(u_j)
+\frac{2(\sum_{j=1}^ku_j-1)_+}{h}\ge1.
\]
\end{lemma}

\begin{proof}We first reduce to the cube $[0,1]^k$.  If one coordinate is
below $0$, increasing it to $0$ cannot increase the displayed expression: on
$(-\infty,0)$ the slope of $G$ is $-2/h$, while the positive-part term
has coordinate slope either $0$ or $2/h$.  If one coordinate is above $1$,
decreasing it to $1$ cannot increase the expression, since both slopes are
nonnegative.  Thus it suffices to work on $[0,1]^k$.

On $[0,1]$ a direct computation gives
\[
G''(s)=-\frac{8h}{(h+2s)^3}<0.
\]
Thus $G$ is concave, and therefore
\[
G(s)\ge(1-s)G(0)+sG(1)
=\frac{2(1-s)}h+\frac{2s}{h+2}.
\]
Writing $U=\sum_j u_j$, we consequently obtain
\[
\sum_{j=1}^kG(u_j)+\frac{2(U-1)_+}{h}
\ge L(U):=\frac{2(k-U)}h+\frac{2U}{h+2}
 +\frac{2(U-1)_+}{h}.
\]
The function $L$ is decreasing on $[0,1]$ and increasing on $[1,k]$;
therefore its minimum is attained at $U=1$. By
\eqref{eq:pd-h-relation}, $2(k-1)/h=h/(h+2)$, and hence
\[
L(1)=\frac{2(k-1)}h+\frac2{h+2}=1.\qedhere
\]
\end{proof}

\begin{proof}[Proof of Theorem~\ref{thm:pd-forward-all-ranks}]
The upper estimate is Proposition~\ref{prop:Construction with better bound}.
For the lower estimate, suppose that $t<\Delta$. By
Lemmas~\ref{lem:pd-forward-endpoint} and \ref{lem:pd-forward-residual},
\[
1=\sum_{j=1}^kA_j+R
>
\sum_{j=1}^kG(u_j)
+\frac{2(\sum_{j=1}^ku_j-1)_+}{h}
\ge1,
\]
where the last inequality follows from
Lemma~\ref{lem:pd-forward-numerical}. This contradiction proves
$t\ge\Delta$.
\end{proof}

\subsection{An endpoint criterion for countable compacta}
\label{subsec:pd-opposite-common}

The following criterion is the only result of this subsection used later.
Since it will be applied only to countable compact spaces, we use the
identification $M(K)=\ell_1(K)$ recorded above.

\begin{prop}[Abstract lower bound for countable compacta]
\label{prop:pd-abstract-opposite-lower}
Let $K,L$ be countable compact spaces, let
$T:\C(L)\to\C(K)$ be a positive unital isomorphism, and put $S=T^{-1}$.
Suppose that there are $e\in L$ and distinct $d_1,\ldots,d_k\in K$, where
$k\ge2$, such that, writing
\[
\Lambda=S^*\delta_e,
\qquad
a_j=(T^*\delta_{d_j})(e),
\]
we have
\begin{equation}\label{eq:pd-opposite-moving-estimate}
a_j>0,
\qquad
\norm{\delta_{d_j}-a_j\Lambda}\le a_j\norm S-1
\quad(1\le j\le k).
\end{equation}
Then $\norm S\ge2+\sqrt3$ if $k=2$, and
$\norm S\ge2+\sqrt5$ if $k\ge3$.
\end{prop}

\begin{proof}
For $x\in K$, put $\mu_x=T^*\delta_x$, and set
$q=\norm\Lambda$, $p_j=\Lambda(d_j)$, and
$r_j=\norm{\delta_{d_j}-a_j\Lambda}$. By unitality, the measures $\mu_x$
are probabilities, $T^*\Lambda=\delta_e$, and $1\le q\le\norm S$. Since
$\mu_{d_j}-a_j\delta_e\ge0$ has norm $1-a_j$,
\begin{equation}\label{eq:pd-opposite-residual-estimate}
r_j=\norm{S^*(\mu_{d_j}-a_j\delta_e)}
\le\norm S(1-a_j).
\end{equation}

We first claim that
\begin{equation}\label{eq:pd-opposite-budgets}
\sum_{p_j>0}p_j(1-a_j)\le\frac{q-1}{2},
\qquad
\sum_{p_j>0}p_j\le\frac{q+1}{2}.
\end{equation}
Indeed, write $\Lambda=\sum_{x\in K}\lambda_x\delta_x$ and put
$c_x=\mu_x(e)$. Since 
\[\sum_x\lambda_x=\Lambda(K)=1 = \delta_e(e) = T^*\Lambda(e) = \sum_x\lambda_xc_x,
\]
we have $\sum_x\lambda_x(1-c_x)=0$. Hence, using $0\leq c_x\leq 1$, $p_j=\lambda_{d_j}$, $a_j=c_{d_j}$, and since the
$d_j$'s are distinct, we obtain
\[
\sum_{p_j>0}p_j(1-a_j)\leq \sum_{\lambda_x>0}\lambda_x(1-c_x)
=\sum_{\lambda_x<0}(-\lambda_x)(1-c_x)
\le \sum_{\lambda_x<0} (-\lambda_x)=:N
\]
and $\sum_{p_j>0}p_j\leq \sum_{\lambda_x>0}\lambda_x=:P$. Further, since $1=\sum_x\lambda_x = P - N$ and $q=\sum_x |\lambda_x| = P + N$, we have $N = \frac{q-1}{2}$ and $P=\frac{q+1}{2}$. Putting everything together, this finishes the proof that \eqref{eq:pd-opposite-budgets} holds.

The identification $M(K)=\ell_1(K)$ also gives
\begin{equation}\label{eq:pd-opposite-coordinate-norm}
r_j=\abs{1-a_jp_j}+a_j(q-\abs{p_j}).
\end{equation}
If $p_j\le0$, then $r_j=1+a_jq\ge1+a_j$, so by \eqref{eq:pd-opposite-moving-estimate} and
\eqref{eq:pd-opposite-residual-estimate} we obtain $a_j\norm S - 1\ge 1+a_j$ and $\norm S(1-a_j)\ge 1+a_j$. This gives $a_j<1$ and
\[
\norm S\ge\max\left\{
\frac{1+a_j}{1-a_j},\frac2{a_j}+1
\right\} =
\max\left\{
\frac2{1-a_j}-1,\,\frac2{a_j}+1
\right\}\ge2+\sqrt5,
\]
where the last inequality follows from Lemma~\ref{lemma:the constant} applied to $(1-a_j,a_j)\in\Delta_2$. Thus all $p_j$ are positive
whenever $\norm S<2+\sqrt5$.

The remainder of the proof uses the following purely numerical claim.

\begin{claim}
Fix $j$ with $p_j>0$, and let $D\ge q$ satisfy
$r_j\le D(1-a_j)$ and $r_j\le a_jD-1$. If $3\le D\le5$, then
\begin{equation}\label{eq:firstFromClaimForCK}
p_j\left(1-a_j+\frac1D\right)
\ge1+\frac{(D+1)(q-D)}{4D}.
\end{equation}
If moreover $D\le2+\sqrt5$, then
$p_j(1-a_j)\ge(q-1)/(D+1)$.
\end{claim}

\begin{proof}[Proof of the Claim]
Once \eqref{eq:pd-opposite-coordinate-norm} is available, this is a purely
numerical argument and does not use any further property of $\C(K)$ spaces.
Put $a=a_j$ and $p=p_j$. By \eqref{eq:pd-opposite-coordinate-norm},
$r_j=\abs{1-ap}+a(q-p)$. Since $\abs{1-ap}\ge1-ap$, we obtain $r_j\ge1+aq-2ap$. Combining this successively with $r_j\le aD-1$ and $r_j\le D(1-a)$ gives
\[
p\ge P_1(a):=\frac{1}{a} + \frac{q-D}{2},
\qquad
p\ge P_2(a):=\frac{q+D}{2} - \frac{D-1}{2a}.
\]
On the other hand, $\abs{1-ap}\ge ap-1$ gives $aq-1\le r_j\le D(1-a)$ which yields $a\le A:=(D+1)/(D+q)\le1$. Further, we have $P_1(a)-P_2(a)=(D+1-2aD)/(2a)$ and so if $P(a)=\max\{P_1(a),P_2(a)\}$, then $P=P_1$ on $(0,a_0]$ and
$P=P_2$ on $[a_0,A]$, where $a_0:=(D+1)/(2D)$ and $a_0\leq A$ using that $q\leq D$.

For $\tau\ge 0$, put $F_\tau(a)=(1-a+\tau)P(a)$. Straightforward computation, using that $a_0\le  (3+1)/6 = 2/3$ since $D\geq 3$, gives for $a\in (0,a_0)$
\[F_\tau'(a)=-\frac{1+\tau}{a^2}+\frac{D-q}{2}\le-\frac{1}{a_0^2}+\frac{D-1}{2}\le -\frac94+\frac{D-1}{2}
\le -\frac94+2
=-\frac14 < 0,\]
while for $a\in (a_0,A)$ we obtain
\[
F_\tau''(a)=-(D-1)(1+\tau)/a^3<0
\]
Consequently,
$F_\tau(a)\ge\min\{F_\tau(a_0),F_\tau(A)\}$.

Set $R=1+(D+1)(q-D)/(4D)$. We compute
\[
a_0P(a_0)
 =
\frac{D+1}{2D}
\left(
\frac{2D}{D+1}+\frac{q-D}{2}
\right) = R,\]
which gives 
\[\begin{split}
F_{1/D}(a_0) & = (1-a_0+\tfrac{1}{D})P(a_0) = a_0P(a_0)=R,\\
F_0(a_0) & = (1-a_0)P(a_0) = \frac{R(1-a_0)}{a_0} = \frac{R(D-1)}{D+1}.
\end{split}\]
Further, we have
\[
P(A)  = \frac{D+q}{2}
  -\frac{(D-1)(D+q)}{2(D+1)} = \frac{D+q}{D+1},
\]
which implies
\[\begin{split}
F_{1/D}(A) & = (1-A+\tfrac{1}{D})P(A) = \frac{q}{D},\\
F_0(A) & = (1-A)P(A) = \frac{q-1}{D+1}.
\end{split}\]

Now 
\[
F_{1/D}(A) - F_{1/D}(a_0) = \frac{q}{D} - R = \frac{(D-3)(D-q)}{4D}\ge0 
\]
and hence
\[
p(1-a+1/D)\ge F_{1/D}(a)\ge \min\{F_{1/D}(A),F_{1/D}(a_0)\} = F_{1/D}(a_0) = R,\]
proving
\eqref{eq:firstFromClaimForCK}.

Finally, if $D\le2+\sqrt5$, then we have $4D-D^2+1\ge 0$ and therefore
\[
F_0(a_0)-F_0(A) = \frac{R(D-1)}{D+1} - \frac{q-1}{D+1}
=\frac{(D-q)(4D-D^2+1)}{4D(D+1)}\ge0.
\]
Therefore $p(1-a)\ge F_0(a)\ge F_0(A)=(q-1)/(D+1)$.
\end{proof}

Suppose first that $k=2$ and $\norm S<2+\sqrt3$, and put
$D=\max\{3,\norm S\}$. Then $3\le D<2+\sqrt3$, and the assumptions of
the Claim are satisfied. Since both $p_j$ are positive,
\eqref{eq:pd-opposite-budgets} and \eqref{eq:firstFromClaimForCK} give
\[
\frac{q-1}{2}+\frac{q+1}{2D}\stackrel{\eqref{eq:pd-opposite-budgets}}{\ge} \sum_{j=1}^2p_j(1-a_j)
  +\frac1D\sum_{j=1}^2p_j  
\stackrel{\eqref{eq:firstFromClaimForCK}}{\ge}2+\frac{(D+1)(q-D)}{2D}.
\]
Hence,
\[
0\le
\frac{q-1}{2}+\frac{q+1}{2D}
-\left(2+\frac{(D+1)(q-D)}{2D}\right)
=\frac{D^2-4D+1}{2D},
\]
a contradiction with the fact that for $D<2+\sqrt{3}$ we have $D^2-4D+1<0$.

Finally, suppose that $k\ge3$ and $\norm S<2+\sqrt5$, and put again
$D=\max\{3,\norm S\}$. Then $3\le D<2+\sqrt5$, so the second part of
the Claim gives
$p_j(1-a_j)\ge(q-1)/(D+1)$ for every $j$. If $q=1$, then \eqref{eq:pd-opposite-budgets} implies $a_j=1$ for every
$j$. Hence \eqref{eq:pd-opposite-residual-estimate} gives
$\delta_{d_j}=\Lambda$ for every $j$, contradicting the fact that the
$d_j$'s are distinct. Thus $q>1$. Summing the inequalities above and
using \eqref{eq:pd-opposite-budgets}, we obtain
\[
\frac{k(q-1)}{D+1}
\le\sum_{j=1}^kp_j(1-a_j)
\le\frac{q-1}{2}.
\]
Consequently, $D\ge2k-1\ge5$, contradicting $D<2+\sqrt5$.
\end{proof}

\begin{remark}
We note that Proposition~\ref{prop:pd-abstract-opposite-lower} remains
valid for arbitrary compact spaces $K$ and $L$, although the proof is
slightly more involved. Since all subsequent applications concern
countable compacta, we have chosen to present only this simpler setting.
\end{remark}

\subsection{The opposite orientation at principal ranks}
\label{sec:pd-higher-opposite}

We now assume that $0<\alpha<\omega_1$ and $\alpha=\omega^\beta$ for some
ordinal $\beta$. Put
\[
K=[1,\omega^\alpha\cdot k],\qquad
L=[1,\omega^\alpha],\qquad
d_j=\omega^\alpha\cdot j\quad(1\le j\le k),\qquad e=\omega^\alpha.
\]
Let $T:\C(L)\to\C(K)$ be a positive unital isomorphism and put
$S=T^{-1}$, $D=\norm S$, $\mu_x=T^*\delta_x$ for $x\in K$,
$\Lambda=S^*\delta_e$, and $a_j=\mu_{d_j}(e)$ for $1\le j\le k$.

The following lemma isolates the part of the higher-rank proof that converts
localized inverse peaks into the moving estimate.
\begin{lemma}[Localized peaks imply the moving estimate]
\label{lem:pd-localized-peaks}
In the notation above, fix $1\le j\le k$. Suppose there are decreasing
clopen neighbourhood bases $(U_n)$ at $e$ in $L$ and $(V_n)$ at $d_j$
in $K$, numbers $\ep_n\downarrow0$, points $z_n\in V_n$, and functions
$f_n\in\C(K,[0,1])$ such that
\[
f_n(z_n)=1,\qquad
\spt f_n\subset V_n,\qquad
\norm{Sf_n|_{L\setminus U_n}}<\ep_n,
\]
and
\[
\sup_{z\in V_n}
\abs{\mu_z(U_n)-\mu_{d_j}(U_n)}<\ep_n.
\]
Then
\[
\frac1D\le a_j\le1,
\qquad
\norm{\delta_{d_j}-a_j\Lambda}\le a_jD-1.
\]
\end{lemma}

\begin{proof}
Write $d=d_j$ and $a=a_j$, and put
$h_n=Sf_n$ and $r_n=\max_{x\in U_n}h_n(x)$. Since $\mu_{z_n}$ is a
probability measure,
\[
1=f_n(z_n)=Th_n(z_n)
 = \int_L h_n d\mu_{z_n}\le r_n\mu_{z_n}(U_n)+\ep_n.
\]
Moreover, $\mu_{z_n}(U_n)\to a$, because
$\mu_d(U_n)\to\mu_d(\{e\})=a$ and
$\abs{\mu_{z_n}(U_n)-\mu_d(U_n)}<\ep_n$. Since $r_n\leq D$ we obtain also $1\leq D\cdot \mu_{z_n}(U_n) + \varepsilon_n\to Da$, so $a\ge 1/D > 0$ and $\liminf_n r_n\ge1/a$. The inequality
$a\le1$ follows from the fact that $\mu_d$ is a probability measure.

Choose $x_n\in U_n$ such that $h_n(x_n)=r_n$, and put
$\Gamma_n=S^*\delta_{x_n}$. Since $x_n\to e$, we have
$\Gamma_n\xrightarrow{w^*}\Lambda$, while $\norm{\Gamma_n}\le D$.
Furthermore, $r_n=\Gamma_n(f_n)\le\Gamma_n^+(V_n)$, because $0\le f_n\le1$ and $\spt f_n\subset V_n$. Hence
Lemma~\ref{lem:pd-shrinking-escape} gives
\[
D\ge\liminf_n\norm{\Gamma_n}
\ge\frac1a+\norm{\Lambda-\frac1a\delta_d}.
\]
Multiplying by $a$ yields
$\norm{\delta_d-a\Lambda}\le aD-1$.
\end{proof}

\begin{lemma}[Principal-rank moving estimate]
\label{lem:pd-higher-endpoint-estimates}
For every $1\le j\le k$,
\[
\frac1D\le a_j\le1,
\qquad
\norm{\delta_{d_j}-a_j\Lambda}\le a_jD-1.
\]
\end{lemma}

\begin{proof}
Fix $j$, and let $(U_n)$ and $(V_n)$ be decreasing clopen neighbourhood
bases at $e$ and $d_j$, respectively. Set $F_n=L\setminus U_n$. Then
$\height(F_n)\le\alpha=\omega^\beta$, whereas
$\height(V_n)=\alpha+1>\omega^\beta$.

Choose $\ep_n\downarrow0$. Since
$z\mapsto\mu_z(U_n)=T\chi_{U_n}(z)$ is continuous, we may shrink $V_n$ so
that
\[
\sup_{z\in V_n}|\mu_z(U_n)-\mu_{d_j}(U_n)|<\ep_n.
\]
Proposition~\ref{prop:infinite splitting old}, applied to
$S:\C(K)\to\C(L)$ with $L_1=V_n$, $U=V_n$, and $L_2=F_n$, provides
$z_n\in V_n$ and $f_n\in\C(K,[0,1])$ such that
\[
f_n(z_n)=1,\qquad
\spt f_n\subset V_n,\qquad
\norm{Sf_n|_{F_n}}<\ep_n.
\]
Lemma~\ref{lem:pd-localized-peaks} now gives the result.
\end{proof}

\begin{thm}\label{thm:pd-higher-opposite-lower}
Let $0<\alpha<\omega_1$ be of the form $\alpha=\omega^\beta$. Then
\[
d_{BM}^+\bigl(\C(\omega^\alpha),\C(\omega^\alpha\cdot2)\bigr)
\ge2+\sqrt3,
\]
and, for every integer $k\ge3$,
\[
d_{BM}^+\bigl(\C(\omega^\alpha),\C(\omega^\alpha\cdot k)\bigr)
\ge2+\sqrt5.
\]
\end{thm}

\begin{proof}
By the normalization at the beginning of the section, it suffices to
consider a positive unital isomorphism $T$ as above. Then
Lemma~\ref{lem:pd-higher-endpoint-estimates} verifies the hypotheses of
Proposition~\ref{prop:pd-abstract-opposite-lower}.
\end{proof}

\begin{proof}[Proof of Theorem~\ref{thm:Intro3}]
For every nonzero countable ordinal $\alpha$ and every $k\ge2$,
Theorem~\ref{thm:pd-forward-all-ranks} gives
\[
d_{BM}^+\bigl(\C(\omega^\alpha\cdot k),\C(\omega^\alpha)\bigr)
=k+\sqrt{(k-1)(k+3)}.
\]

Suppose now that $\alpha=\omega^\beta$ for some countable ordinal $\beta$.
Theorem~\ref{thm:pd-higher-opposite-lower} gives the lower bounds
\[
d_{BM}^+\bigl(\C(\omega^\alpha),\C(\omega^\alpha\cdot2)\bigr)
\ge2+\sqrt3
\]
and, for $k\ge3$,
\[
d_{BM}^+\bigl(\C(\omega^\alpha),\C(\omega^\alpha\cdot k)\bigr)
\ge2+\sqrt5.
\]
For $k=2$, the matching upper bound follows from
Proposition~\ref{prop:pd-higher-two-upper}. For $k\ge3$,
Proposition~\ref{prop:constructing-isomorphism-new-attempt}(i) applied
with $n=1$ and $\beta = k$, gives the
matching upper bound $2+\sqrt5$.
\end{proof}

\begin{remark}
The assumption $\alpha=\omega^\beta$ is used only in
Lemma~\ref{lem:pd-higher-endpoint-estimates}, where we apply
Proposition~\ref{prop:infinite splitting old} with
$\omega^\beta=\alpha$. For arbitrary $\alpha$, we would need the same
proposition with an arbitrary ordinal in place of $\omega^\beta$, but this
is false already for the ordinal $2$. Indeed, take
$K_1=L_1=U=[1,\omega^2]$ and $K_2=L_2=[1,\omega]$. Then
$\height(L_1)=3>2=\height(L_2)$, while $\C(K_1)$ and $\C(K_2)$ are
isomorphic, since both are isomorphic to $c_0$. If
$R:\C(K_1)\to\C(K_2)$ is an isomorphism, the claimed extension would give,
for every $\ep>0$, a norm-one $f\in\C(K_1)$ such that
$\norm{Rf}<\ep$, contradicting $\norm{Rf}\ge\norm{R^{-1}}^{-1}$.
Thus the required localization result does not extend to arbitrary
ordinals, although the corresponding extension of
Theorem~\ref{thm:pd-higher-opposite-lower} may still be true.
\end{remark}

\section{Concluding remarks and open problems}

Theorem~\ref{thm:Intro3} settles the forward positive distances at every
countable rank. In the opposite orientation, our argument gives the exact
values at rank one and at principal higher ranks, leaving the following
problem.

\begin{question}\label{que:pd-opposite-arbitrary-rank}
Let $0<\alpha<\omega_1$. Is it true that
\[
d_{BM}^+\bigl(\C(\omega^\alpha),\C(\omega^\alpha\cdot2)\bigr)
=2+\sqrt3
\]
and, for every $k\ge3$,
\[
d_{BM}^+\bigl(\C(\omega^\alpha),\C(\omega^\alpha\cdot k)\bigr)
=2+\sqrt5?
\]
\end{question}

In Question~\ref{qe:milj} we asked about a generalization of Miljutin's
theorem. All proofs of Miljutin's theorem known to us rely on the
Pe{\l}czy\'nski decomposition method, which does not appear to extend to
positive isomorphisms. Two particularly important instances of
Question~\ref{qe:milj} are the following.

\begin{question}

\begin{itemize}
    \item Is there a positive isomorphism of $\C([0, 1])$ onto $\C(2^{\omega})$?
    \item Is there a positive isomorphism of $\C(2^{\omega})$ onto $\C([0, 1])$?
\end{itemize}
\end{question}

We note that for the non-surjective case, the answer is positive as it is known that an isomorphic embedding of $\C(2^{\omega})$ in any Banach lattice $E$ implies the existence of a positive isomorphic embedding (\cite[Theorem A]{Ghoussoub_positive_embeddings_C(Delta)}), and linear embeddability of $\C([0, 1])$ in $E$ is equivalent to even lattice embeddability (\cite[Theorem A]{Aviles-Cervantes-RuedaZoca-Tradacete_linearvslatticeembeddings}). 

We also consider the following nonlinear version of our result. We recall that a long-standing open problem in the nonlinear theory of Banach spaces asks whether the Szlenk index is invariant under bi-Lipschitz bijections of Banach spaces. An important result hinting at a possible positive answer is \cite{Dutrieux}.
The problem remains open also when restricted to just the class of $\C(K)$ spaces, and even when the considered compact spaces are countable. It is true when one of the spaces is isomorphic to $c_0$ (see \cite[Theorem 2.2]{Godefroy2000SubspacesO}), in this case even for uniform homeomorphisms instead of bi-Lipschitz bijections (see \cite[Theorem 5.5(i)]{GodefroyKaltonLancien_uniform}).
In spite of these important steps forward and a substantial effort of the mathematical community, the problem remains wide open. We suggest considering the assumption that the bi-Lipschitz bijection is order-preserving (which for linear mappings coincides with being positive), which might shed some new light on the question. 

\begin{question}
Do order-preserving bi-Lipschitz bijections of $\C(K)$ spaces over countable compact spaces $K$ preserve the Szlenk index of $\C(K)$? Do they preserve even the height of $K$?
\end{question}

Related to this question, let us recall Proposition~\ref{prop:cAndc0} by which there does not exist an order-preserving bijection from $c\simeq \C([1, \omega])$ to $c_0\simeq \C_0([1, \omega])$ with Lipschitz inverse. On the other hand, $c_0$ is positively isomorphic to $c$. This is witnessed by the linear operator $T:\C_0([1, \omega]) \rightarrow \C([1, \omega])$, defined by
\[
Tf(\omega)=f(1), \quad Tf(n)=f(1)+f(n+1), \quad n \in \en,
\]
which is a positive surjective isomorphism with inverse $S:\C([1, \omega]) \rightarrow \C_0([1, \omega])$,
\[
Sh(1)=h(\omega), \quad Sh(n+1)=h(n)-h(\omega), \quad n \in \en,
\qquad Sh(\omega)=0.
\]

Another question related to our research concerns spaces of continuous functions with the pointwise topology. We recall that the classification of spaces of the form $\C_p(K)$, where $K$ is a countable compact space, by linear homeomorphism, coincides with the classification of the Banach spaces $\C(K)$ by linear isomorphisms, see \cite{BG88}. Thus, the following problem seems natural.

\begin{question}
Assume that $K$ and $L$ are countable compact spaces. Is it true that there
exists a positive linear homeomorphism of $\C_p(K)$ onto $\C_p(L)$ if and
only if there exists a positive isomorphism of $\C(K)$ onto $\C(L)$?
\end{question}

% Bibliography entries whose metadata was corrected in this revision.
\makeatletter
\newcommand{\bluebibentry}[1]{\expandafter\gdef\csname bluebib@#1\endcsname{1}}

\let\originalbibitem\bibitem
\renewcommand{\bibitem}[1]{%
  \par\color{black}\originalbibitem{#1}%
  \@ifundefined{bluebib@#1}{}{\color{blue}}%
}
\makeatother

\let\L\PolishL

\section*{Statements and Declarations}

\noindent\textbf{Funding.}
M. C\'uth was supported by the Czech Science Foundation, project no. GA\v{C}R 24-10705S.

\medskip
\noindent\textbf{Competing interests.}
The authors have no relevant financial or non-financial interests to disclose.

\medskip
\noindent\textbf{Data availability.}
No datasets were generated during the current study.

\medskip
\noindent\textbf{Use of generative artificial intelligence.}
The main ideas underlying the lower-bound arguments in
Section~\ref{sec:positive-distances}, which yield the exact positive
Banach--Mazur distances stated in Theorem~\ref{thm:Intro3}, emerged
through an extended interaction with ChatGPT (models 5.5 and
5.6). The tool was used to explore proof strategies, perform and check
calculations, identify gaps in preliminary arguments, and refine the
final proofs. All arguments were critically reviewed and independently
verified by the authors, who take full responsibility for the
mathematical content.
\bibliography{references.bib}\bibliographystyle{siam}
\end{document}